\documentclass[openany, a4paper, 12pt]{article}

\usepackage{amsfonts}

 \usepackage{mathrsfs,amsfonts,amsmath}
 \usepackage{color}

 \makeatletter
 \@addtoreset{equation}{section}
 \makeatother
 \newtheorem{theorem}{Theorem}[section]
 \newtheorem{definition}{Definition}[section]
 \newtheorem{hypothesis}{Hypothesis}[section]
 \newtheorem{lemma}{Lemma}[section]
 \newtheorem{proposition}{Proposition}[section]
 \newtheorem{corollary}{Corollary}[section]
 \newtheorem{remark}{Remark}[section]
 \newtheorem{example}{Example}[section]

 \def\beqlb{\begin{eqnarray}}\def\eeqlb{\end{eqnarray}}
 \def\beqnn{\begin{eqnarray*}}\def\eeqnn{\end{eqnarray*}}

 \def\mbb{\mathbb}

 \def\qed{\hfill$\Box$\medskip}

\newcommand{\bcen}{\begin{center}}
\newcommand{\ecen}{\end{center}}
\newcommand{\bgeqn}{\begin{equation}}
\newcommand{\edeqn}{\end{equation}}

\def\dz{\delta}

\def\wz{\wedge}
\def\ez{\epsilon}

\def\lra{\longrightarrow}

\def\rar{\rightarrow}

\def\ra{\rangle}
\def\la{\langle}

\begin{document}
\begin{center}{\Large\textbf{Rescaled Lotka-Volterra Models
Converge to Super Stable Processes}\footnote{Supported by NSFC
(No.10721091 )}}\end{center}

\bigskip

\centerline{ Hui He\footnote{ \textit{E-mail address:} {
hehui@bnu.edu.cn }} }

\medskip

\centerline{Laboratory of Mathematics and Complex Systems, }

\centerline{ School of Mathematical Sciences, Beijing Normal
University,}

\smallskip

\centerline{ Beijing 100875, People's Republic of China}

\bigskip

{\narrower{\narrower{\narrower

\begin{abstract}
Recently, it has been shown that stochastic spatial Lotka-Volterra
models when suitably rescaled can converge to a super Brownian
motion. We show that the limit process could be a super stable
process if the kernel of the underlying motion is in the domain of
attraction of a stable law. The corresponding results in Brownian
setting were proved by Cox and Perkins (2005, 2008). As applications
of the convergence theorems, some new results on the asymptotics of
the voter model started from single 1 at the origin are obtained
which improve the results by Bramson and Griffeath (1980).
\end{abstract}

\smallskip

\noindent\textit{AMS 2000 subject classifications.} Primary 60K35,
60G57; secondary 60F17, 60J80.

\smallskip

\noindent\textit{Key words and phrases.} {super stable process},
{Lotka-Volterra}, {voter model}, {domain of attraction}, {stable
law}, {stable random walk}.

\smallskip

\noindent\textbf{Abbreviated Title:} Lotka-Volterra model

\par}\par}\par}

\bigskip\bigskip

\section{Introduction}

 \subsection{Motivation}\label{secmotivation} Originally, super Brownian
motion arises as the limit of branching random walks; see
\cite{[D93], [CDP01], [P02]}. Recently, it has been shown that many
interacting particle systems with very different dynamics, when
suitably rescaled, all converge to super Brownian motion. Such
examples include  the voter model, the contact process, interacting
diffusion process and  the spatial Lotka-Volterra model; see
\cite{[CDP01],[DP99],[CK03], [CP05],[CP08]}. Donsker's invariance
principle is deeply involved in those results; see \cite{[S02]} for
an excellent nontechnical introduction. So if we assume that the
kernel of the underlying motion has finite variance, super Brownian
motion is obtained as the limit process. On the other hand, the
general class of stable distribution was introduced and given this
name by the famous French mathematician Paul L\'{e}vy. The
inspiration for L\'{e}vy was the desire to generalize the Central
Limit Theorem which is the foundation of Donsker's principle. Thus
we can expect that if we let the kernel of the underlying motion be
in the domain of attraction of a stable law, the limit process could
be a super stable process.

 A motivation for proving those limit theorems is to actually
use it in the study of complicated approximating systems.
 For example, the Lotka-Volterra
invariance principle established in \cite{[CP05]} was used  to study
the coexistence and survival problem of the Lotka-Volterra model;
see \cite{[CP07]}. Cox and Perkins \cite{[CP04]} used the voter
invariance principle to give a probabilistic proof of the
asymptotics for the voter model obtained in \cite{[BG80]}.  In this
paper, we will show that rescaled stochastic spatial Lotka-Volterra
models can converge to super stable processes and also use those
limit theorems to get some new results on the asymptotics for the
voter model.  Coexistence and survival for the Lotka-Volterra model
will be discussed in a future work.

\subsection{Our model}\label{ourmodel}

A stochastic spatial version of the Lotka-Volterra model was first
introduced and studied by Neuhauser and Pacala \cite{[NP99]}. In
this paper, we follow the construction of the model suggested by
\cite{[CP05]} but we assume that the kernel of the model is in the
domain of attraction of a symmetric stable law. We first briefly
describe the model. Let $\{p(x,y)\}$ be a random walk kernel on
${\mbb{Z}}^d$ (the $d$-dimensional integer lattice). Suppose at each
site of ${\mbb Z}^d$ there is a plant of one of two type. We label
the two types 0 and 1. At random times plants die and are replaced
by new plants. The times and the types depend on the configuration
of surrounding plants. We denote by $\xi_t$, an element of
$\{0,1\}^{{\mbb Z}^d}$,  the state of the system at time $t$ and
$\xi_t(x)$ gives the type of the plant at $x$ at time $t$. To
describe the evolution of the system, for $\xi\in\{0,1\}^{{\mbb
Z}^d}$, define
 \bgeqn
 \label{1.1}
f_i(x,\xi)=\sum_{y\in {\bf Z}^d}p(x,y)1_{\{\xi(y)=i\}},\quad i=0,1.
 \edeqn
 Let $\alpha_0$, $\alpha_1$ be nonnegative
parameters. Define the Lotka-Volterra \textit{rate function}
$c(x,\xi)$ by
\begin{eqnarray*}
c(x,\xi)=\left\{\begin{array}{lll}f_1(f_0+\alpha_0f_1)\quad\textrm{if
}\xi(x)=0,\\
f_0(f_1+\alpha_1f_0)\quad\textrm{if }\xi(x)=1. \end{array}\right.
\end{eqnarray*}
The Lotka-Volterra process $\xi_t$ is the unique
$\xi\in\{0,1\}^{{\bf Z}^d}$-valued Feller process with rate function
$c(x,\xi)$, meaning that the generator of $\xi_t$ is the closure of
the operator $\Omega$
$$
\Omega\phi(\xi)=\sum_xc(x,\xi)(\phi(\xi^x)-\phi(\xi))
$$
on the set of function $\phi:\xi\in\{0,1\}^{{\mbb
Z}^d}\rightarrow\mbb R$ depending on only finitely many coordinates,
where $\xi^x(y)=\xi(y)$ for $y\neq x$ and $\xi^x(x)=1-\xi(x)$.

Note that $f_0+f_1=1$. The dynamics of $\xi_t$ can now be described
as follows: at site $x$ in configuration $\xi$, the coordinate
$\xi(x)$ makes transitions
\begin{eqnarray*}
0\rightarrow1\quad\quad\textrm{at rate}\quad
f_1(f_0+\alpha_0f_1)=f_1+(\alpha_0-1)f_1^2,\cr
1\rightarrow0\quad\quad\textrm{at rate}\quad
f_0(f_1+\alpha_1f_0)=f_0+(\alpha_1-1)f_0^2.
\end{eqnarray*}
These rates are interpreted in \cite{[NP99]} as follows. A plant of
type $i$ t site $x$ dies at rate $f_i+\alpha_if_{1-i}$, and is
replaced by a plant of type $\xi(y)$ where $y$ is chosen with
probability $p(x,y)$. $\alpha_i$ measures the strength of
interspecific competition of type $i$ and we set the
self-competition parameter equal to one.

In \cite{[CDP01]} an invariance principle was proved for the voter
model. That is appropriately rescaled voter models converge to
super-Brownian motion. Thus we can expect that when the parameters
$\alpha_i$ are close to one a similar result holds for the
Lotka-Volterra model. The results in \cite{[CP05]} and \cite{[CP08]}
say that it is true. The intuition of the voter invariance principle
is that when appropriately rescaled, the dependence on the local
density of particles gets washed out and the rescaled voter models
should behave like the rescaled branching random walk. The
asymptotics behavior of the latter is well known: it approaches
super-Brownian motion. On the other hand, if the kernel of the
underlying  motion is in the domain of attraction of a stable law,
appropriately rescaled branching random walk could approach a super
stable  process; see Theorem II.5.1 of \cite{[P02]}. The above
reasoning suggests the possibility of that suitably rescaled
Lotka-Volterra should approach a super stable process. Our main
results in this paper will show that it is the case.

Let $M(\mbb R^d)$ denote the space of finite measures on $\mbb R^d$,
endowed with the topology of weak convergence of measures. Let
$\Omega_D=D([0,\infty),M(\mbb R^d))$ be the Skorohod space of
c\`{a}dl\`{a}g paths taking values  in $M(\mbb R^d)$. Let $\Omega_C$
be the space of continuous $M(\mbb R^d)$-valued paths with the
topology of uniform convergence on compact set. We denote by
$X_t(\omega)=\omega_t$ the coordinate function. We write $\mu(\phi)$
for $\int\phi d\mu$. For $1\leq n\leq\infty$ let $C_b^n(\mbb R^d)$
be the space of bounded continuous function whose partial
derivatives of order less than $n+1$ are also bounded and
continuous, and let $C_0^n(\mbb R^d)$ be the space of those
functions in $C_b^n(\mbb R^d)$ with compact support.

A $\mbb R^d$-valued L\'{e}vy process $Y_t$ is said to be a symmetric
$\alpha$-stable process with index $\alpha\in(0,2]$ and diffusion
speed $\sigma^2>0$ if \bgeqn\label{stablelaw}\Psi(\eta):={
E}(e^{i\eta\cdot Y_1})=e^{-\sigma^2|\eta|^\alpha},\edeqn where $|y|$
is the Euclidean norm of $y$. The distribution of $Y_1$ will be
called $(\sigma^2,\alpha)$-stable law. When $\alpha=2$, $Y_t\in \mbb
R^d$ is a $d$-dimensional $\sigma^2$-Brownian motion whose generator
is ${\cal A}\phi=\frac{\sigma^2\Delta\phi}{2}$ for $\phi\in
C_b^2({\mbb R^d})$. When $0<\alpha<2$, the generator of $Y_t$ is
given by
$$
{\cal A}\phi(x)=\frac{\sigma^2\Delta^{\alpha/2}\phi(x)}{2}=
\sigma^2\int\left[\phi(x+y)-\phi(x)-\frac{1}{1+|y|^2}
\sum_{i=1}^dy_jD_j\phi(x)\right]\nu(dy)
$$
for $\phi\in C_b^2(\mbb R^d)$ and $D_j=\frac{\partial}{\partial
x_j}$, where $$\nu(dy)=c|y|^{-d-\alpha}1_{\{|y|\neq0\}}(dy)$$ for an
appropriate $c>0$; see \cite{[S99]} for details. In both cases,
$C_b^{\infty}(\mbb R^d)$ is a core for $\cal A$ in that the
$bp$-closure of $\{(\phi,{\cal A}\phi):\phi\in C_b^{\infty}\}$
contains $\{(\phi,{\cal A}\phi):\phi\in{\cal D}({\cal A})\}$, where
${\cal D}({\cal A})$ denotes the domain of the weak generator for
the process $Y$; see \cite{[P02]}.

An adapted a.s.-continuous $M(\mbb R^d)$-valued process
$\{X_t:t\geq0\}$ on a complete filtered probability space $(\Omega,
{\cal F}, {\cal F}_t, P)$ is said to to a \textit{super symmetric
$\alpha$-stable process with branching rate $b\geq0$, drift
$\theta\in \mbb R$ and diffusion coefficient $\sigma^2>0$ starting
at $X_0\in M(\mbb R^d)$} if it solves the following martingale
problem:

\noindent{\narrower For all $\phi\in C_b^{\infty}(\mbb R^d)$,
\begin{eqnarray} \label{MP1}
M_t(\phi)=X_t(\phi)-X_0(\phi)-\int_0^tX_s\left(
\frac{\sigma^2\Delta^{\alpha/2}\phi(x)}{2} \right)ds
-\theta\int_0^tX_s(\phi)ds
\end{eqnarray}
is a continuous $({\cal F}_t)$-martingale, with $M_0(\phi)=0$ and
predictable square function \bgeqn\label{MP2} \la
M(\phi)\ra_t=\int_0^tX_s(b\phi^2)ds.\edeqn
\par}
\noindent The existence and uniqueness in law of a solution to this
martingale problem is well known; see Theorem II.5.1 and Remark
II.5.13 of \cite{[P02]}. Let $P^{b,\theta,\sigma^2,\alpha}_{X_0}$
denote the law of the solution on $\Omega_C$. So $b$ and  $\theta$
can be regarded as branching parameters and parameters $\sigma$ and
$\alpha$ determine the underlying motion.

Let $\{Z_n:n\geq1\}$ be a discrete time random walk on $\mbb{Z}^d$,
$$
Z_n=z_0+\sum_{i=1}^nU_i,
$$
where $z_0\in\mbb{Z}^d$ and the random variables $(U_i:i\geq1)$ are
independent identically distributed on $\mbb{Z}^d$. Let $\{p(x,y)\}$
be a random walk kernel. In the following of this paper we assume
that

\noindent {\narrower {\bf(A1)}: $p(x,y)=p(x-y)$ is an irreducible,
symmetric, random walk kernel on ${\mbb{Z}}^d$ and $p(0)=0$. For
$\alpha\in(0,2]$ and  $\sigma^2>0$, $\{p(x)\}$ is in the domain of
attraction of a symmetric $(\sigma^2,\alpha)$-stable law; i.e.,
$$ P(U_1=x)=p(x)
$$
and there exists a function $b(n)$ of regular variation of index
$1/\alpha$ such that
 \bgeqn \label{defDA}
 b(n)^{-1}\sum_{i=1}^nU_i\xrightarrow{(d)}Y_1\quad\textrm{ as
 }n\rightarrow\infty,
 \edeqn where $Y_1$ is determined by (\ref{stablelaw}) and the symbol
 $\xrightarrow{(d)}$ means convergence in distribution.
\par}
\noindent We will call a random walk (discrete time or continuous
time) with kernel satisfying assumption (A1) a stable random walk.
In the following of this paper, we always assume that
$$\textit{(A1) holds for some}~ \sigma>0 ~and ~\alpha\in(0,2].$$

\begin{remark}
\label{remark1.1} Without loss of generality, we may and will assume
that function $b$ is continuous and monotonically increasing from
$\mbb R^{+}$ onto $\mbb R^{+}$ and $b(0)=0$; see \cite{[LR91]} or
\cite{[F71]}. We also have that
$$
b(x)=x^{1/\alpha}s(x),\quad x>0,
$$
where $s:(0,\infty)\rightarrow(0,\infty)$ is a slowly varying
function, meaning that for any $c>0$,
$$
\lim_{x\rightarrow\infty}\frac{s(cx)}{s(x)}=1
$$
where the convergence holds uniformly when $c$ varies over the
interval $[\ez, 1/\ez]$ for any $\ez>0$; see Lemma 2 of VIII.8 of
\cite{[F71]}.
\end{remark}
\begin{remark}
\label{transient} According to Proposition 2.5 of \cite {[LR91]} and
its proof, we have that under (A1), random walk $\{Z_n\}$ is
transient if and only if
$$
\sum_{k=1}^{\infty}b(k)^{-d}<\infty.
$$
By Lemma 2 in Section VIII.8 of \cite{[F71]}, the random walk is
always transient when $d>\alpha$.
 Typically, when $d=\alpha=1$,
the random walk is recurrent if only if  $$
\sum_{k=1}^{\infty}\frac{1}{ks(k)}=\infty. $$
\end{remark}

 Now, we are ready to define our rescaled
Lotka-Volterra models. For $N=1,2,\cdots,$ let
$$\mbb{S}_N=\mbb{Z}^d/b(N).$$ Define the kernel $p_N$ on
$\mbb{S}_N$ by
$$
p_N(x)=p(xb(N)),\quad\quad x\in\mbb{S}_N.
$$
For $\xi\in\{0,1\}^{\mbb{S}_N}$, define the densities
$f_i^N=f_i^N(\xi)=f_i^N(x,\xi)$ by
$$
f_i^N(x,\xi)=\sum_{y\in\mbb{S}_N}p_N(y-x)1_{\{\xi(y)=i\}},\quad\quad
i=0,1.
$$
Let $\alpha_i=\alpha_i^N$ depend on $N$ and let $\xi_t^N$ be the
process taking values in $\{0,1\}^{\mbb{S}_N}$ determined by the
rates: at site $x$ in configuration $\xi$, the coordinate $\xi(x)$
makes transitions
\begin{eqnarray*}
0\rightarrow1\quad\quad\textrm{at rate}\quad
Nf_1^N(f_0^N+\alpha_0^Nf_1^N),\cr 1\rightarrow0\quad\quad\textrm{at
rate}\quad Nf_0^N(f_1^N+\alpha_1^Nf_0^N).
\end{eqnarray*}
That is $\xi^N_t$ is rate-$N$ Lotka-Volterra process determined by
the parameters $\alpha_i^N$ and the kernel $p_N$. More precisely, if
set
\begin{eqnarray*}
c_N(x,\xi)=\left\{\begin{array}{lll}Nf_1^N(f_0^N+\alpha_0^Nf_1^N)
\quad\textrm{if
}\xi(x)=0,\\
Nf_0^N(f_1^N+\alpha_1^Nf_0^N)\quad\textrm{if }\xi(x)=1,
\end{array}\right.
\end{eqnarray*}
$\xi^N_t$ is the unique Feller process taking values in
$\{0,1\}^{\mbb S^N}$ whose generator is the closure of the operator
$$
\Omega_N\phi(\xi)=\sum_{x\in\mbb
S^N}c_N(x,\xi)(\phi(\xi^x)-\phi(\xi))
$$
on the set of function $\phi:\xi\in\{0,1\}^{{\mbb
Z}^d}\rightarrow\mbb R$ depending on only finitely many coordinates.
Here $\xi^x(y)=\xi(y)$ for $y\neq x$ and $\xi^x(x)=1-\xi(x)$.

 \begin{remark} If we assume
$\sum_{x\in\mbb{Z}^d}x^ix^jp(x)=\dz_{ij}\sigma^2<\infty$, then
$p(x)$ is in the domain of attraction of a normal law. That is
 the case of $\alpha=2$. So we recover the fixed kernel models in
\cite{[CP05]}. For critical case, since there are significant
differences between the case of $d=\alpha=1$ and the case of
$d=\alpha=2$,  we only consider the case of $d=\alpha=1$. For
$d=\alpha=2$, please see the work in \cite{[CP08]}. \end{remark}

Define $$g(x)=\int_1^x b(s)^{-1}ds$$ for $d=\alpha=1$ and $x\geq0$.
According to Remark \ref{transient}, the one-dimensional random walk
$Z$ is recurrent if and only if $\lim_{x\rar\infty} g(x)=\infty.$

 Set
\begin{eqnarray*}
N'=\begin{cases}N,&\textrm{if }d>\alpha,\\
N,&\textrm{if }d=\alpha=1 \textrm{ and }
\lim_{x\rar\infty}g(x)<\infty ,\\
N/{g(N)},&\textrm{if }d=\alpha=1\textrm{ and }\lim_{x\rar\infty}
g(x)=\infty.
\end{cases}
\end{eqnarray*}
That is when the stable random walk is transient $N'=N$ and
$N'=N/g(N)$ if the stable random walk is recurrent.

We define the corresponding measure-valued process $X_t^N$ by
 \bgeqn
 \label{resMV}
 X_t^N=\frac{1}{N'}\sum_{x\in\mbb{S}_N}\xi_t^N(x)\dz_x.
 \edeqn
 As in \cite{[CP05]} and \cite{[CP08]}, we make the following
 assumptions:
 \begin{align}
 \label{A2}
  &(1)~\sum_{x\in\mbb{S}_N}\xi_0^N(x)<\infty.\notag\\
&(2)~X_0^N\rightarrow X_0\quad\quad \textrm{in }~ M(\mbb R^d)
 \quad\textrm{ as }N\rightarrow\infty. \tag{${\bf A2}$}\\
 &(3)~\theta_i^N=N'(\alpha_i^N-1)\rightarrow\theta_i\in\mbb R\quad\quad
 \textrm{as }N\rightarrow\infty,\quad i=0,1.\notag
 \end{align}
Now, we are ready to describe our main results.
\subsection{Main results}\label{MR} To describe the limit process,
 we introduce a coalescing
random walk systems $\{\hat{B}_t^x,x\in \mbb{Z}^d\}$. Each
$\hat{B}_t^x$ is a rate 1 random walk on $\mbb{Z}^d$ with kernel
$p$, with $\hat{B}_0^x=x$. The walks move independently until they
collide, and then move together after that. For finite
$A\subset\mbb{Z}^d$, let $$\hat{\tau}(A)=\inf\{s:|\{\hat{B}_t^x,x\in
A\}|=1\}$$ be the time at which the particles starting from $A$
coalesce into a single particle, and write $\hat{\tau}(a,b,\cdots)$
when $A=\{a,b,\cdots\}$. Note that when the stable random walk is
transient, we can define the ``escape'' probability by
$$
\gamma_e=\sum_{e\in \mbb{Z}^d}p(e)P(\hat{\tau}(0,e)=\infty).
$$
We also define
 \begin{eqnarray*}
 &&\beta=\sum_{e,e'\in\mbb Z^d}p(e)p(e')
 P(\hat{\tau}(e,e')<\infty,\hat{\tau}(0,e)=\hat{\tau}(0,e')=\infty),\cr
 &&\dz=\sum_{e,e'\in\mbb Z^d}p(e)p(e')
 P(\hat{\tau}(0,e)=\hat{\tau}(0,e')=\infty).
 \end{eqnarray*}
We also need a collection of independent (noncoalescing) rate-1
continuous time random walks with step function $p$, which we will
denote $\{B_t^{x}:x\in \mbb Z^d \}$, such that $B_0^{x}=x$. Define
the collision times $$ \tau(x,y)=\inf\{t\geq0:B_t^x=B_t^y\},\quad
x,y\in \mbb Z^d.$$ Let $P_N$ denote the law of $X^N_.$. Our first
result is following.
\begin{theorem}
\label{mainUP} Assume (A1), (A2) and $d\geq\alpha$. If the stable
random walk is transient,  then
$$P_N\xrightarrow{(d)}P_{X_0}^{2\gamma_e,\theta,\sigma^2,\alpha}$$ as
$N\rightarrow\infty$, where $\theta=\theta_0\beta-\theta_1\dz$.
\end{theorem}

Note that if we assume
$\sum_{x\in\mbb{Z}^d}x^ix^jp(x)=\dz_{ij}\sigma^2<\infty$, then
$\{p(x)\}$ is in the domain of attraction of a normal law with
$b(N)=\sqrt{N}$.  So Theorem \ref{mainUP} generalizes Theorem 1.2 in
\cite{[CP05]}.

Next, we consider the recurrent case.  And for some technical
reasons we need to assume that the $\{p(x)\}$ is in the \textit{
domain of normal attraction} of $(\sigma^2,1)$-stable law; see
Remark \ref{LVreasoning} below. To state our result, we introduce
the one-dimensional potential kernel $a(x)$,
 \bgeqn
 \label{Pa}
 a(x)=\int_0^{\infty}\left[P(B_t^0=0)-P(B_t^x=0)\right]dt.
 \edeqn
We will discuss the existence of $a(x)$ later.  Note that
$a(x)\geq0$. Let $\{p_t(x):t\geq0,x\in\mbb R\}$ denote the
transition density of  $\{Y_t\}$. Now we define
 \begin{eqnarray}
 \label{Gamma}
 \gamma^{\ast}=(p_1(0))^{-1}\int_0^\infty\sum_{x,y,e,e'}p(e)p(e')
 P({\tau}(0,e)&\wz&{\tau}(0,e')>{\tau}(e,e')\in du,\cr
 && B_u^0=x,B_u^e=y)a(y-x).
 \end{eqnarray}
 Our critical Lotka-Volterra invariance principle is
 \begin{theorem}
 \label{main2}
 Assume  (A2), $d=\alpha=1$,  (A1) holds with $b(t)=t$ and $N'=N/\log
 N$. Then
 $$P_N\xrightarrow{(d)}P_{X_0}^{2\hat{p},\theta,\sigma^2,1}$$ as
$N\rightarrow\infty$, where
$\theta=\gamma^{\ast}(\theta_0-\theta_1)$ and
$\hat{p}=(p_1(0))^{-1}$.
 \end{theorem}

\begin{remark}
According to Remark \ref{transient}, the assumption that (A1) holds
with $b(t)=t$ implies that the stable random walk is recurrent.
\end{remark}


Now, we consider the applications of the convergence theorems. One
can see from the rate function form that if we set
$\alpha_0=\alpha_1=1$, $\xi_t$ is just the  well known \textit{voter
model}. Identify $\xi_t$ with the set $\{x:\xi_t(x)=1\}$ and let
$\xi_t^A$ denote the voter model starting from 1's exactly on $A$,
$\xi_0^A=A$. Write $\xi_t^x$ for $\xi_t^{\{x\}}$. The usual additive
construction of the voter models yields
$$
\xi_t^A=\bigcup_{x\in A}\xi_t^x.
$$
The fact that $|\xi_t^0|=\sum_x\xi_t^0(x)$ is martingale tells us
$|\xi_t^0|$ hits 0 eventually with probability 1. Letting
$p_t=P(|\xi_t^0|>0)$, it follows that $p_t\rar0$ as $t\rar\infty$.
People always want to determine the rate at which $p_t\rar0$. By
using a result in \cite{[Sa79]}, Bramson and Griffeath \cite{[BG80]}
were able to obtain precise asymptotics under the assumption that
the underlying motion is a simple random walk. By making the voter
model invariance principle, Cox and Perkins \cite{[CP04]} reproved
the main result in \cite{[BG80]} under a weaker assumption that the
jump kernel has finite variance. In this paper as applications of
the convergence theorems above we want to determine the rate at
which $p_t\rar0$ under the assumption (A1). With notation $f(t)\sim
g(t)$ as $t\rar\infty$ we mean $\lim_{n\rar\infty}f(t)/g(t)=1$. Our
result is following theorem.
 \begin{theorem}
 \label{voterasy}
 Assume $d\geq\alpha$ and (A1) holds with $b(t)=t^{1/\alpha}$;
 i.e., $\{p(x)\}$ is in
 the domain of normal attraction of the $(\sigma,\alpha)$-stable
 law. Let $\gamma_{1}=p_1(0)^{-1}$ for $d=\alpha$.
 Then as $t\rar\infty$
 \begin{alignat*}{2}
 p_t&\sim\frac{\log t}{\gamma_1t}\qquad &d=\alpha,\\
    &\sim{(\gamma_et)^{-1}}\qquad &d>\alpha.
 \end{alignat*}
Moreover,
 $$
 P\left(p_t|\xi_t^0|>u\big{|}|\xi_t^0|>0\right)\xrightarrow{t\rar\infty}
  e^{-u},\quad u>0.
 $$
 \end{theorem}

\medskip

At last, we introduce some notations which will play important roles
in our proofs of the main results. First, according to \cite{[F71]},
for $0<\underline{\alpha}\leq\alpha$, we can define
$$
|p|_{\underline{\alpha}}:=\sum_{x\in\mbb
Z^d}|x|^{\underline{\alpha}}p(x)<\infty.
$$ And by (A2), define
$$\bar{\theta}=1\vee\sup_{N,i}N'|\alpha_i^N-1|<\infty.$$
For $D\subset\mbb R^d$ and $\phi: D\rightarrow\mbb R$, define
$$
||\phi||_{\textrm{Lip}}=||\phi||_{\infty}+\sup_{x\neq
y}\frac{|\phi(x)-\phi(y)|}{|x-y|}.
$$
For $0<\underline{\alpha}\leq1$, let \beqnn
||\phi||_{\underline{\alpha}}=\left\{\begin{array}{ccc}0,&\phi\equiv
c\text{ for some constant } c\in\mbb R\\
\sup_{x\neq
y,|x-y|\leq1}\frac{|\phi(x)-\phi(y)|}{|x-y|^{\underline{\alpha}}}
\vee2||\phi||_{\infty},& \text{ otherwise},\end{array}\right. \eeqnn
and for $\underline{\alpha}>1$ let
$$||\phi||_{\underline{\alpha}}=2||\phi||_{\text{Lip}}.$$
Note that for $\underline{\alpha}\leq 1$,
 $$
\sup_{x\neq
y,|x-y|\leq1}\frac{|\phi(x)-\phi(y)|}{|x-y|^{\underline{\alpha}}}
\leq\sup_{x\neq y}\frac{|\phi(x)-\phi(y)|}{|x-y|}.
 $$
Thus for any $\underline{\alpha}>0$ { \bgeqn
 \label{LIP}
||\phi||_{\underline{\alpha}}\leq
2||\phi||_{\textrm{Lip}}\quad\textrm{and }\quad|\phi(x)-\phi(y)|\leq
||\phi||_{\underline{\alpha}}|x-y|^{\underline{\alpha}}.
 \edeqn}
\begin{remark}
\label{holder} Since $p(\cdot)$ in this paper may not have bounded
moment of the first order, we can not use Lipschitz norm to do
estimates. Thus a `H\"{o}lder' norm is introduced.
\end{remark}

\medskip

The remaining of this paper is organized  as follows. In Section
\ref{Pre}, we first give some random walk estimates and then deduce
the semimartingale decompositions for the approximating processes.
Finally, we prove a key result, uniform convergence of random walk
generators to the generator of the symmetric stable process. In
Section 3 and Section 4, we follow the strategy in \cite{[CP05]} and
\cite{[CP08]} to prove our convergence theorems, Theorem
\ref{mainUP} and Theorem \ref{main2}. Our proofs will be deeply
involved due to the lack of high moments. We will carry out in
detail only the part that differs. Theorem \ref{voterasy} will be
proved in Section 5.

\section{Preliminaries}\label{Pre}
\subsection{Random walk estimate}
Recall that $\{B_t^x,x\in\mbb Z^d\}$ is a collection of rate-one
independent stable random walks with $B_0^x=x$. Let
$p_t(x,y)=P(B_t^x=y)$ denote the transition function of $\{B_t^x\}$.
We denote by $l$ the inverse of $b$. Define the characteristic
function of the step function $p(\cdot)$ by
$$
\psi(\eta)=\sum_xp(x)e^{-iy\cdot\eta}\quad\textrm{ for }\quad\eta\in
T^d:=(-\pi,\pi]^d.
$$
Since $p$ is symmetric, $\psi(\eta)$ is real. So
 \beqlb\label{A3}
 p_t(0,x)\leq p_t(0,0).
 \eeqlb
The following proposition is taken from \cite{[LR91]}.
 \begin{proposition}
 \label{DAeqv}
 The following are equivalent:
 \begin{enumerate}
 \item[(1)]$p(\cdot)$ is in the domain of attraction of
 $(\sigma^2,\alpha)$-stable law.

 \item[(2)]$\psi(\eta)=1-\frac{\sigma^2}{l(1/|\eta|)}
 +o\left(\frac{1}{l(1/|\eta|)}\right)$
 as $|\eta|$ tends to 0.

 \item[(3)]
 $\psi\left(\frac{\eta}{b(n)}\right)^n\xrightarrow{n\rightarrow\infty}\Psi(\eta)$,
 \quad$\eta\in\mbb R^d$.
\end{enumerate}
\end{proposition}

We also have that $l$ is of regular variation of index $\alpha$ and
$$
l(x)=x^{\alpha}t(x),
$$
where
$$
t(x)=s(l(x))^{-\alpha}.
$$
By Lemma 2.1 in \cite{[LR91]}, for any $\ez>0$, we have that there
exist two positive constants $C_{\ez}, C'_{\ez}$ such that, for any
$1\leq y\leq z$,
 \bgeqn
 \label{regular}
C_{\ez}y^{\alpha-\ez}\leq l(y)\leq C'_{\ez}y^{\alpha+\ez} \quad
\textrm{ and }\quad
C_{\ez}\left(\dfrac{z}{y}\right)^{\alpha-\ez}\leq\frac{l(z)}{l(y)}\leq
C_{\ez}'\left(\dfrac{z}{y}\right)^{\alpha+\ez}.
 \edeqn
A similar result also holds for $b$, with $\alpha$ replaced by
$1/\alpha$. Since $p(\cdot)$ is symmetric and irreducible, $\psi$ is
real and $\psi(\eta)=1$ if and only if $\eta=0$; see \cite{[Sp76]}.
According to Proposition \ref{DAeqv}, we may assume that there
exists a constant $C>0$ such that
 $$
\frac{C}{l(1/|\eta|)}\leq1-\psi(\eta)\leq1
 $$
for every $\eta\in T^d$. (\ref{regular}) tells us  that for
$b(t)\geq d\pi$, and $0\leq\ez\leq\alpha,$
 \bgeqn \label{boundCha}
 t(1-\psi(\frac{\eta}{b(t)}))\geq \frac{Cl\left(b(t)\right)}{l\left(b(t)/|\eta|\right)}\geq
(C_{\ez}\vee C_{\ez}')(|\eta|^{\alpha+\ez}+|\eta|^{\alpha-\ez}).
 \edeqn
Recall that $\{p_t(x):t\geq0,x\in\mbb R\}$ denote the transition
density of $\{Y_t\}$.  The local limit theorem for the stable random
walk which plays an important role in our proofs of main results
will be given in the following proposition.
 \begin{proposition}
 \label{CPtostable}
 If (A1) holds,
 \bgeqn
 \label{tranapp}
 \lim_{t\rightarrow\infty}\sup_{x\in\mbb Z^d}
 \left|b(t)^dp_t(0,x)-p_1\left(\frac{x}{b(t)}\right) \right|=0
 \edeqn
and there exists a constant $C$ depending on $p(\cdot)$ such that
for every $t\geq0$ and $x\in \mbb
 R^d$,
 \bgeqn
 \label{boundtransition}
 p_t(0,x)\leq C b(t)^{-d}.
 \edeqn
 Moreover,  if $b(t)=t$ and $d=1$,
 \bgeqn
 \label{7.6}
 \sup_{x\in \mbb Z}P(B_t^0=x)\leq C_{\ref{7.6}}(t+1)^{-1}.
 \edeqn
 \end{proposition}
{\bf Proof.} Since $l$ is a function of regular variation,
 by Proposition \ref{DAeqv}, for each $|\eta|>0$,
 \begin{eqnarray}
 \label{concha}
 \lim_{t\rightarrow\infty}
  t\left(1-\psi\left(\frac{\eta}{b(t)}\right)\right)
  =\lim_{t\rightarrow\infty}\frac{l(b(t))}{l(b(t)/|\eta|)}
 (\sigma^2+o(1))=\sigma^2|\eta|^{\alpha}.
 \end{eqnarray}
Then
 \begin{eqnarray*}
 &&\left|b(t)^dp_t(0,x)-p_1\left(\frac{x}{b(t)}\right) \right|\\
  &&\quad\leq(2\pi)^{-d}\left|\int_{b(t)T^d}e^{-ix\cdot(\eta/b(t))}
  \exp\left\{-t\left(1-\psi\left(\frac{\eta}{b(t)}\right)\right)\right\}d\eta
   -\int_{b(t)T^d}e^{-i(x/b(t))\cdot\eta}\Psi(\eta)d\eta\right|\\
   &&\quad\quad+(2\pi)^{-d}\int_{\mbb R^d\setminus
   b(t)T^d}\exp\left\{-\sigma^2|\eta|^{\alpha}\right\}d\eta\\
   &&\quad\leq(2\pi)^{-d}\int_{b(t)T^d}\left|
  \exp\left\{-t\left(1-\psi\left(\frac{\eta}{b(t)}\right)\right)\right\}
   -\exp\left\{-\sigma^2|\eta|^{\alpha}\right\}\right|d\eta\\
   &&\quad\quad+(2\pi)^{-d}\int_{\mbb R^d\setminus
   b(t)T^d}\exp\left\{-\sigma^2|\eta|^{\alpha}\right\}d\eta.
 \end{eqnarray*}
Then the Dominated Convergence Theorem with (\ref{boundCha}) yields
(\ref{tranapp}). For (\ref{boundtransition}), when $b(t)\geq d\pi$,
 \begin{eqnarray*}
 p_t(0,x)&=&(2\pi)^{-d}\int_{T^d}e^{-ix\cdot\eta}
  \exp\left\{-t\left(1-\psi(\eta)\right)\right\}d\eta\\
 &\leq&(2\pi)^{-d}b(t)^{-d}\int_{b(t)T^d}
  \exp\left\{-t\left(1-\psi\left(\frac{\eta}{b(t)}\right)\right)\right\}d\eta\\
  &\leq&(2\pi)^{-d}b(t)^{-d}\int_{\mbb R^d}
  \exp\{-(C_{\ez}\vee
  C_{\ez}')(|\eta|^{\alpha+\ez}+|\eta|^{\alpha-\ez})\}d\eta\\&\leq&
  Cb(t)^{-d},
 \end{eqnarray*}
where the second inequality follows from (\ref{boundCha}). Then
(\ref{boundtransition}) holds for every $t\geq0$. We complete the
proof. \qed

\smallskip

 The following two propositions consider the growth of the
stable random walk.
\begin{proposition}
\label{lemma7.3} (a) If $z_T\in \mbb Z^d$ and $t_T>0$ satisfy
 \bgeqn
 \label{7.8}
 \lim_{T\rar\infty}\frac{z_T}{b(T)}=z \textrm{ and }\lim_{T\rar\infty}
 \frac{t_T}{T}=s>0
 \edeqn
 then
 \bgeqn
 \label{7.9}
 \lim_{T\rar\infty}b(T)^dP(B_{t_T}^0=z_T)=\frac{p_1(z/s)}{s^d}.
 \edeqn

 (b) For each $K>0$, there is a constant $C_{\ref{7.10}}(K)>0$ such
 that
 \bgeqn
 \label{7.10}
 \liminf_{T\rar\infty}\inf_{|x|\leq Kb(T)}b(T)^dP(B_T^0=x)\geq
 C_{\ref{7.10}}(K).
 \edeqn
\end{proposition}
{\bf Proof. } By (\ref{7.8}) and Remark \ref{remark1.1}, we have
$\lim_{T\rar\infty}\frac{b(t_T)}{b(T)}=s$. Then (\ref{7.9}) follows
from (\ref{tranapp}).  For (b), when $\alpha=2$, by (\ref{tranapp}),
the desired result is immediate. When $0<\alpha<2$, recall that
$\{p_t(x):t\geq0, x\in \mbb R^d\}$ is the transition density of a
symmetric $\alpha$-stable process. By the arguments after Remark 5.3
of \cite{[BL02]}, there exists two positive constants $c_1$ and
$c_2$ such that
 \bgeqn
 \label{estistable}
c_1\left(t^{-d/\alpha}\wedge\frac{t}{|x|^{d+\alpha}}\right)\leq
p_t(x)\leq
c_2\left(t^{-d/\alpha}\wedge\frac{t}{|x|^{d+\alpha}}\right).
 \edeqn
By above bounds and (\ref{tranapp}),
 \beqnn
\liminf_{T\rar\infty}\inf_{|x|\leq Kb(T)}b(T)^dP(B_T^0=x)& =&
\liminf_{T\rar\infty}\inf_{|x|\leq Kb(T)}p_1(x/b(T))\\
&\geq& c\left(1\wedge K^{d+\alpha}\right).
 \eeqnn
 The desired result
follows readily. \qed

\smallskip

\begin{proposition}
\label{propincrease} Assume $d=1$. If $g_1$ and $g_2$ are two
positive functions on $\mbb R^+$ such that
$g_1(x)\rar+\infty,~g_2(x)\rar+\infty$ as $x\rar+\infty$, then there
is exists a constant $C_{\ref{increasein}}$ which only depends on
$p$ such that
 \beqlb\label{increasein}
 P\left(|B_{g_1(N)}^0|\geq g_2(N)\right)\leq
 \frac{C_{\ref{increasein}}g_1(N)}{l(g_2(N))}.
 \eeqlb
\end{proposition}
{\bf Proof. }First,
$$P\left(|B_{g_1(N)}^{0}|\geq g_2(N)\right)\leq
P\left(\max_{ u\leq {g_1(N)}}|B_u^0|\geq g_2(N)\right).$$  Note that
$\{B_u^0:u\geq0\}$ is a compound Poisson process whose L\'{e}vy
measure is given by
$$
 \nu_0(dz):=\sum_{y\in \mbb Z^d}p(y)\dz_y(dz),
$$
which is a symmetric measure. According to the arguments in Section
3 of \cite{[Pr81]},
$$
P\left(\max_{ u\leq {g_1(N)}}|B_u^0|\geq g_2(N)\right)\leq
Cg_1(N)\left(\nu_0(z:|z|>g_2(N))+g_2(N)^{-2}\int_{|z|\leq
{g_2(N)}}z^2\nu_0(dz)\right),
$$where $C$ is a positive constant; see (3.2) of \cite{[Pr81]}.
Since $p(\cdot)$ is in the domain of attraction of
$(\sigma,\alpha)$-stable law, we have
 \bgeqn
 \label{for3}
\frac{x^2[\nu_0(z:|z|>x)]}{\int_{|z|\leq
x}z^2\nu_0(dz)}\lra\frac{2-\alpha}{\alpha}
 \edeqn
and
 \bgeqn
 \label{for2}
\frac{x\int_{|z|\leq b(x)}z^2\nu_0(dz)}{b(x)^2}\lra C_0
 \edeqn
as $x\rar\infty$ for some constant $C_0>0$; see (5.16) and (5.23) in
Chapter XVII of \cite{[F71]}. By (\ref{for3}) there exists a
constant $C_1$ independent  $N$ such that
$$
\nu_0\left(z:|z|>g_2(N)\right)\leq C_1{g_2(N)^{-2}}{\int_{|z|\leq
g_2(N)}z^2\nu_0(dz)}.
$$
According to (\ref{for2}), there exists another constant $C_2$
independent of $N$ such that
 \beqnn
 {g_2(N)^{-2}}{\int_{|z|\leq g_2(N)}z^2\nu_0(dz)}
 \leq
 \frac{C_2}{l(g_2(N))}.
 \eeqnn
(Recall that $l$ is the inverse function of $b$.) Thus
$$
P\left(\max_{ u\leq {g_1(N)}}|B_u^0|\geq g_2(N)\right)\leq
CC_2(C_1+1)\frac{g_1(N)}{l(g_2(N))}
$$
which yields the desired result. \qed

\subsection{Semimartingale decompositions}\label{secsemi}
Some results in this subsection are exactly the same with those in
Section 3 of \cite{[CP08]}. For complement, we list them here. Let
$\xi_t^N$ be the rescaled Lotka-Volterra model we have constructed
in Section \ref{ourmodel}. As in \cite{[CP08]}, we introduce the
following notation. If
$$
\phi=\phi_s(x),\quad \dot{\phi}_s(x)\equiv\frac{\partial}{\partial
s}\phi(s,x)\in C_b([0,T]\times\mbb S_N),
$$
and $s\leq T$, define
 \begin{eqnarray}
 \label{generatorN}
 {\cal A}_N(\phi_s)(x)&=&\sum_{y\in \mbb
 S_N}Np_N(y-x)(\phi_s(y)-\phi_s(x))\\
 D_t^{N,\,1}(\phi)&=&\int_0^tX_s^N({\cal
 A}_N\phi_s+\dot{\phi_s})ds\\
 D_t^{N,\,2}(\phi)&=&\frac{N(\alpha_0^N-1)}{N'}\int_0^t
 \sum_{x\in\mbb
 S_N}\phi_s(x)1_{\{\xi_s^N(x)=0\}}(f_1^N(x,\xi_s^N))^2ds\\
 D_t^{N,\,3}(\phi)&=&\frac{N(\alpha^N_1-1)}{N'}\int_0^t
 \sum_{x\in\mbb
 S_N}\phi_s(x)1_{\{\xi_s^N(x)=1\}}(f_0^N(x,\xi_s^N))^2ds
 \end{eqnarray}
 \begin{eqnarray}
 \la M^N(\phi)\ra_{1,\,t}&=&\frac{N}{(N')^2}\int_0^t\sum_{x\in \mbb
 S_N}\phi_s^2(x)\sum_{y\in\mbb
 S_N}p_N(y-x)(\xi_s^N(y)-\xi_s^N(x))^2ds\\
\la M^N(\phi)\ra_{2,\,t}&=&\frac{1}{(N')^2}\int_0^t\sum_{x\in \mbb
 S_N}\phi_s^2(x)\big{[}(\alpha_0^N-1)1_{\{\xi_s^N(x)=0\}}
 (f_1^N(x,\xi_s^N))^2\cr
 &&\quad
 +(\alpha_1^N-1)1_{\{\xi_s^N(x)=1\}}(f_0^N(x,\xi_s^N))^2\big{]}ds
 \end{eqnarray}
If $X_{\cdot}$ is a process let $({\cal F}_t^X,t\geq0)$ be the
right-continuous filtration generated by $X_{\cdot}$. The following
proposition is a version of Proposition 3.1 of \cite{[CP08]}. For
its proof, please go to Section 2 of \cite{[CP05]}.
 \begin{proposition}
 \label{semidec}
 For $\phi, \dot{\phi}\in C_b([0,T]\times\mbb S_N)$ and $t\in
 [0,T]$,
 \bgeqn
 \label{XNdec}
 X_t^N(\phi_t)=X_0^N(\phi_0)+D_t^N(\phi)+M_t^N(\phi),
 \edeqn
 where
 \bgeqn
 \label{Driftdec}
 D_t^N(\phi)=D_t^{N,1}(\phi)+D_t^{N,2}(\phi)-D_t^{N,3}(\phi)
 \edeqn
and $M_t^N(\phi)$ is an ${\cal F}_t^{X^N}$-square-integrable
martingale with predictable square function
 \bgeqn
 \label{premartN}
 \la M^N(\phi)\ra_t=\la M^N(\phi)\ra_{1,t}+\la M^N(\phi)\ra_{2,t}.
 \edeqn
 \end{proposition}

The following lemma is a generalization of Lemma 3.5 of
\cite{[CP05]} and Lemma 4.8 of \cite{[CP08]}.
\begin{lemma}
\label{estimateMNt} There is a constant $C$ such that if
$\phi:[0,T]\times\mbb S_N\rightarrow\mbb R$ is a bounded measurable
function, then

(a) $\la M^N(\phi)\ra_{2,t}=\int_0^tm_{2,s}^N(\phi)ds$, where
 \bgeqn
 \label{m2st}
 |m_{2,s}^N(\phi)|\leq C\frac{||\phi_s||_{\infty}^2}{(N')^2}X_s^N(1).
 \edeqn

(b) For $\underline{\alpha}<1\wedge\alpha,$ \bgeqn \label{MN1tdec}
    \la
    M^N(\phi)\ra_{1,t}=2\int_0^tX_s^N((N/N')\phi_s^2f_0^N(\xi_s^N))ds
    +\int_0^tm_{1,s}^N(\phi_s)ds,
    \edeqn
where
 \bgeqn
 \label{m1sNest}
 |m_{1,s}^N(\phi)|\leq
 \left[X_s^N(1)\frac{2N||\phi||^2_{\underline{\alpha}}
 |p|_{\underline{\alpha}}}{N'b(N)^{\underline{\alpha}}}
 \right]
 \wedge\left[\frac{2N||\phi||^2_{\infty}X_s^N(1)}{N'}\right].
 \edeqn

(c) For $i=2,3, D_t^{N,i}(\phi)=\int_0^td_s^{N,i}(\phi)ds$ for
$t\leq T$, where for all $N$, $s\leq T$,
 $$
|d_s^{N,i}(\phi)|\leq C||\phi_s||_{\infty}X_s^N\left(
(N/N')f_0^N(\xi_s^N)\right).
 $$
\end{lemma}
\begin{remark}
Note that when $N'=N$, since $f_0^N\leq1$,
$$
|d_s^{N,i}(\phi)|\leq C||\phi_s||_{\infty}X_s^N( 1),\quad i=2,3.
$$
\end{remark}
{\bf Proof.} (a)  In the following of this proof, with $C$ we denote
a positive constant which may change from line to line. Since
$f_0^N\leq1$, $f_1^N\leq1$ and $1_{\{\xi_s^N(x)=1\}}=\xi_s^N(x)$,
the definition of $\la M^N(\phi)\ra_{2,t}$ and the fact that
$f_0^N+f_1^N=1$ imply
 \beqnn
 |m_{2,s}^N(\phi)|&\leq&\frac{||\phi||^2_{\infty}\sup_N N'|\alpha_0^N-1|}
 {(N')^3} \sum_{x\in\mbb
 S_N}(f_1^N(x,\xi_s^N))1_{\{\xi_s^N(x)=0\}}
 \\&&
 +\frac{||\phi||^2_{\infty}\sup_N N'|\alpha_1^N-1|}{(N')^2}X_s^N(1)\\
 &\leq&\frac{C||\phi||^2_{\infty}}{(N')^3}\sum_{x,y}p_N(x-y)
 (1-1_{\{\xi_s^N(x)=1\}})1_{\{\xi_s^N(y)=1\}}
 +\frac{C||\phi||^2_{\infty}}{(N')^2}X_s^N(1)\\
 &\leq&\frac{C||\phi||^2_{\infty}}{(N')^2}X_s^N(1),
 \eeqnn
where the second inequality follows from (A2).
 For (b), note that
 \beqnn
 \la M^N(\phi)\ra_{2,t}&=&\frac{1}{(N')^2}\int_0^t\sum_{x\in \mbb
 S_N}\phi_s^2(x)\sum_{y\in\mbb
 S_N}Np_N(y-x)(\xi_s^N(y)-\xi_s^N(x))^2ds\\
&=&\frac{1}{(N')^2}\int_0^t\sum_{x\in \mbb
 S_N}\phi_s^2(x)\sum_{y\in\mbb
 S_N}Np_N(y-x)\left(2\xi_s^N(x)(1-\xi_s^N(y))\right)ds\\
 &&\quad+\frac{1}{(N')^2}\int_0^t\sum_{x\in \mbb
 S_N}\phi_s^2(x)\sum_{y\in\mbb
 S_N}Np_N(y-x)\left(\xi_s^N(y)-\xi_s^N(x)\right)ds.
 \eeqnn
Thus (\ref{MN1tdec}) holds with
 \beqnn
m_{1,s}^N(\phi)&=&\frac{N}{(N')^2}\sum_{x\in\mbb
S_N}\phi^2_s(x)\sum_{y\in \mbb S_N}p_N(y-x)(\xi_s^N(y)-\xi_s^N(x))\\
&=&\frac{N}{(N')^2}\sum_{x\in\mbb S_N}\phi^2_s(x)\sum_{y\in \mbb
S_N}p_N(y-x)(\xi_s^N(y)1_{\{\xi_s^N(x)=0\}}-\xi_s^N(x)1_{\{\xi_s^N(y)=0\}})\\
&=&\frac{N}{(N')^2}\sum_{x,y\in\mbb S_N}
p_N(y-x)(\phi^2_s(x)-\phi_s^2(y))\xi_s^N(y)(1-\xi_s^N(x))\\
&\leq&\frac{2N||\phi||^2_{\infty}X_s^N(1)}{N'}.
 \eeqnn
On the other hand, $$|\phi^2_s(x)-\phi_s^2(y)|\leq
2||\phi||^2_{\underline{\alpha}}|x-y|^{\underline{\alpha}}$$ for
$\underline{\alpha}<1\wedge\alpha.$ Thus
 \beqnn
 m_{1,s}^N(\phi)&\leq&
2(N/N')||\phi||^2_{\underline{\alpha}}\frac{1}{N'}
\sum_y\xi_s^N(y)\sum_x|y-x|^{\underline{\alpha}}p_N(y-x)\\
&\leq&
 X_s^N(1)\frac{2N||\phi||^2_{\underline{\alpha}}
 |p|_{\underline{\alpha}}}{N'b(N)^{\underline{\alpha}}}.
 \eeqnn
 We complete the proof of (b).
For (c), according to (A2), the fact that both $f_0^N$ and $f_1^N$
are less than 1 yields
 \beqnn
 |d_s^{N,i}(\phi)|&\leq&\frac{N\sup_N N'|\alpha^N_{i-2}-1|}{N'}||\phi_s||_{\infty}
 \frac{1}{N'}\sum_x\sum_yp_N(y-x)\xi_s^N(x)(1-\xi_s^N(y))\\
 &\leq& C||\phi_s||_{\infty}X_s^N((N/N')f_0^N(\xi_s^N)).
  \eeqnn
We are done.
 \qed

\subsection{Convergence of Generators}
 In this subsection we consider the uniform convergence
of ${\cal A}_N$. Recall the definition of generators of symmetric
stable processes and the stable random walk $Z_n$ defined in section
\ref{ourmodel}. For each $N>1$, let $\{P_t^{(N)}:t\geq0\}$ be a
rate-$N$ Poisson process which is independent of $\{U_i:i\geq1\}$.
Then
$$
\hat{Z}_t^N=b(N)^{-1}\sum_{i=1}^{P_t^{(N)}}U_i
$$
is a compound Poisson process on $\mbb R^d$ whose L\'{e}vy measure
is given by
 $$
 \nu_N(dy):=\sum_{z\in \mbb S_N}Np_N(z)\dz_z(dy);
 $$
 see \cite{[S99]}.
 Note that both the law of $\hat{Z}_1^N$
and the $(\sigma^2,\alpha)$-stable law are infinitely divisible
distributions. We also have that
$$
{\bf E}\left(e^{-i\hat{Z}_1^N\cdot\eta}\right)
=\exp\left\{-N\left(\psi\left(\frac{\eta}{b(N)}\right)-1\right)\right\}.
$$
By (\ref{concha}),
$$
\hat{Z}_1^N\xrightarrow{(d)}Y_1\quad\textrm{as}\quad N\rar\infty.
$$
According to Theorem 8.7 of \cite{[S99]} and its proof, we see
$$
\rho_N(dy):=\frac{|y|^2}{1+|y|^2}\nu_N(dy)\rightarrow\rho(dy):=
\frac{\sigma^2|y|^2}{1+|y|^2}\nu(dy)\quad\textrm{ in }M(\mbb R^d).$$
For $f\in C_b(\mbb R^d)$, define
$$
||f||_{BL}=\sup_x|f(x)|\vee\sup_{x\neq y}\frac{|f(x)-f(y)|}{|x-y|}.
$$
Let  ${\cal P},{\cal Q}$ be two probability measures on $\mbb R^d$.
Set
$$
||{\cal P}-{\cal Q}||_{BL}:=\sup_{||f||_{BL}=1}\left|\int fd{\cal
P}-\int fd{\cal Q}\right|.
$$
It is easy to see that
 \bgeqn
 \label{probcon2}
||{\cal P}-{\cal Q}||_{BL}=\sup_{||f||_{BL}<\infty}\frac{\left|\int
fd{\cal P}-\int fd{\cal Q}\right|}{||f||_{BL}}.
 \edeqn
 By Problem 3.11.2 of \cite{[EK86]},
 \bgeqn\label{probcon}
 ||{\cal P}-{\cal Q}||_{BL}\leq 3 {\cal M}({\cal P},{\cal Q}),
 \edeqn
where $\cal M$ denotes the Prohorov metric; see Chapter 3 of
\cite{[EK86]}.

\begin{lemma}
\label{conge} For $\phi\in C_b^{1,3}([0,T]\times\mbb R^d),$
$$
\lim_{N\rightarrow\infty}\sup_{s\leq T}||{\cal
A}_N\phi_s-\frac{\sigma^2\Delta^{\alpha/2}\phi_s}{2}||_{\infty}=0.
$$
Moreover, for each $R<\infty$, the rate of convergence is uniform on
 \begin{eqnarray*}
 H_R:=\left\{\phi\in C_b^{1,3}([0,T]\times\mbb R^d):
   \sup_{s,i,j,k}(||\phi_s||_{\infty}+||(\phi_s)_i||_{\infty}
   +||(\phi_s)_{ij}||_\infty+||(\phi_s)_{ijk}||_\infty)<R\right\},
 \end{eqnarray*}
 where the subscripts $i,j,k$ indicate partial derivatives with
 respect to the spatial variable.
\end{lemma}
{\bf Proof.} Recall that  $D_j=\frac{\partial}{\partial x_j}$.
Define
$$
g_s(x,y)=\left[\phi_s(x+y)-\phi_s(x)-\frac{1}{1+|y|^2}
\sum_{i=1}^dy_jD_j\phi_s(x)\right]\cdot \frac{1+|y|^2}{|y|^2}.
$$
 Since $p_N$ is
symmetric, we may rewrite
$$
{\cal A}_N\phi_s(x)=\int g_s(x,y)\rho_N(dy)
$$
and we also have that
$$
\frac{\sigma^2\Delta^{\alpha/2}\phi_s(x)}{2}=\int g_s(x,y)\rho(dy).
$$
Let  $h:\mbb R^d\rightarrow[0,1]$ be a $C_b^{\infty}$ function such
that
$$
B(0,1)\subset\{x:h(x)=0\}\subset\{x:h(x)<1\}\subset B(0,2)
$$
and
$$
B(0,2)^c\subset\{x:h(x)=1\}.
$$
Define $ h_k(x)=h(kx)$ for $k\geq1$. Let
$$
g_k(s,x,y):=h_k(y)g_s(x,y).
$$
Then $g_k(s,x,y)=g_s(x,y)$ for $|y|>2/k$. One can check that
$$
\sup_k\sup_{\phi\in H_R}\sup_s\sup_x\left(||g_k(s,x,\cdot||_{\infty}
+||g_s(x,\cdot)||_{\infty}\right)<C_d R
$$
and for each $k\geq1$
$$
\sup_{\phi\in H_R}\sup_s\sup_x||\sum_{j=1}^d |\frac{\partial
g_k(s,x,y)}{\partial y_j}|||_{\infty}<kC_d R,
$$
where $C_d$ is a constant which only depend on $d$. Typically, for
each $k\geq1$,
$$
\sup_{\phi\in H_R}\sup_s\sup_x||g_k(s,x,\cdot)||_{BL}<(k+1)C_d R.
$$
By (\ref{probcon2}) and (\ref{probcon}), we obtain
\begin{eqnarray*}
&&\sup_{\phi\in H_R}\sup_{s\leq T}\sup_x\left|\frac{\int
g_k(s,x,y)\rho_N(dy)}{\rho_N(\mbb R^d)}-\frac{\int
g_k(s,x,y)\rho(dy)}{\rho(\mbb R^d)}\right|\\&&\quad\quad\leq(k+1)C_d
R \cdot3{\cal M}\left(\frac{\rho_N}{\rho_N(\mbb
R^d)},\frac{\rho}{\rho(\mbb
R^d)}\right)\\&&\quad\quad\rightarrow0,\quad \textrm{ as
}N\rightarrow\infty.
\end{eqnarray*}
By triangle inequality,
 \begin{eqnarray*}
&&\sup_{\phi\in H_R}\sup_{s\leq T}\sup_x\left|{\int
g_k(s,x,y)\rho_N(dy)}-{\int
g_k(s,x,y)\rho(dy)}\right|\\&&\quad\quad\leq C_dR\left|\rho_N(\mbb
R^d)-\rho(\mbb R^d)\right|\\&&\quad\quad\quad+\rho(\mbb
R^d)\sup_{\phi\in H_R}\sup_{s\leq T}\sup_x\left|\frac{\int
g_k(s,x,y)\rho_N(dy)}{\rho_N(\mbb R^d)}-\frac{\int
g_k(s,x,y)\rho(dy)}{\rho(\mbb
R^d)}\right|\\&&\quad\quad\rightarrow0,\quad \textrm{ as
}N\rightarrow\infty.
 \end{eqnarray*}
 Using triangle inequality again,
 \begin{eqnarray*}
&& \sup_{\phi\in H_R}\sup_{s\leq T}||{\cal
A}_N\phi_s-\frac{\sigma^2\Delta^{\alpha/2}\phi_s}{2}||_{\infty}\\
&&\quad\quad \leq\sup_{\phi\in H_R}\sup_{s\leq T}\sup_x\left|\int
g_s(x,y)\rho_N(dy)-\int
g_k(s,x,y)\rho_N(dy)\right|\\&&\quad\quad\quad
 +\sup_{\phi\in H_R}\sup_{s\leq
T}\sup_x\left|\int g_k(s,x,y)\rho_N(dy)-\int
g_k(s,x,y)\rho(dy)\right|\\&&\quad\quad\quad
 +\sup_{\phi\in H_R}\sup_{s\leq
T}\sup_x\left|\int g_k(s,x,y)\rho(dy)-\int g_s(x,y)\rho(dy)\right|\\
&&\quad\quad\leq C_d R\rho_N(\{y:|y|\leq2/k\})+C_d
R\rho(\{y:|y|\leq2/k\})\\&&\quad\quad\quad
 +\sup_{\phi\in H_R}\sup_{s\leq
T}\sup_x\left|\int g_k(s,x,y)\rho_N(dy)-\int
g_k(s,x,y)\rho(dy)\right|
 \end{eqnarray*}
Note that $\rho(dy)$ is absolutely continuous with respect to the
Lebesgue measure. Letting $N$ go to infinity above yields
$$
\lim_{N\rightarrow\infty}\sup_{\phi\in H_R}\sup_{s\leq T}||{\cal
A}_N\phi_s-\frac{\sigma^2\Delta^{\alpha/2}\phi_s}{2}||_{\infty}
\leq2C_d R\rho(\{y:|y|\leq2/k\}).
$$
Then since $\rho(\{0\})=0$ the  desired result follows readily if we
let $k\rightarrow\infty$. \qed

\section{Proof of Theorem \ref{mainUP}}

In this section, we assume the stable random walk $Z$ is transient,
which is equivalent to
$$\int_1^{\infty}\frac{dx}{b(x)^d}<\infty.$$
When $d=\alpha=1$, above condition implies that $s(x)\rar\infty$ as
$x\rar\infty$.The strategy of the proof is  the same with that used
in \cite{[CP05]}.
In \cite{[CP05]} the authors worked with a more general class of
particle systems they called voter perturbations. As a result we
will specialize the setting there for the reader's convenience.
 Let $\{\hat{B}_t^{N,x}:x\in\mbb S_N\}$ denote a
rate-$N$ continuous time coalescing random walk system on $\mbb S_N$
with step function $p_N$ such that $\hat{B}_0^{N,x}=x$. For a finite
set $A\subset \mbb S_N$, let
$$\hat{\tau}^N(A)=\inf\{t\geq0:|\{\hat{B}_t^{N,x},x\in A\}|=1\}.$$
We also need a collection of independent (noncoalescing) rate-$N$
continuous time random walks on $\mbb S_N$ with step function $p_N$,
which we will denote $\{B_t^{N,x}:x\in \mbb S_N \}$, such that
$B_0^{N,x}=x$.
 For any finite subset $A$ of $\mbb
Z^d$, let $\hat{\tau}^N(A)=\hat{\tau}(A/b(N))$. We first check the
kernel assumptions in Section 1.2 of \cite{[CP05]}.
 \begin{lemma}
 \label{Kerlemma}
 There exists a positive sequence $\{\ez_N^{*}\}$ with
 $\ez_N^*\rightarrow0$ and $N\ez_N^*\rightarrow\infty$. such that
 the following hold:
 \begin{eqnarray}
 \lim_{N\rar\infty}NP(B_{\ez_N^*}^{N,0}=0)&=&0.\\
 \lim_{N\rar\infty}\sum_{e\in\mbb
 S_N}p_N(e)P(\hat{\tau}^N(\{0,e\})\in(\ez_N^*,t])&=&0\quad \textrm{ for all
 }\qquad t>0,\cr
 \lim_{N\rar\infty}\sum_{e\in\mbb
 S_N}p_N(e)P(\hat{\tau}^N(\{0,e\})>\ez_N^*)&=&\gamma_e.
 \end{eqnarray}
and if we define $\sigma_N(A)=P(\hat{\tau}^N(A)\leq\ez_N^*)$ for any
finite subset $A$ of $\mbb Z^d$, then
\bgeqn\label{K3}\lim_{N\rar\infty}\sigma_N(A)=\sigma(A)\quad\textrm{
exists}. \edeqn
 \end{lemma}
{\bf Proof }. First, consider the case $d>\alpha$. We may assume
$\ez_N^*=N^{-\ez^*}$ for some $0<\ez^*<1$. We need to find a
suitable condition on $\ez^*$. Recall that $b$ is a function of
regular variation with index $1/\alpha$. Given $\ez<1/2$, there
exist two positive constants $C_{\ez}$, $C_{\ez}'$ such that for
$y\geq1$,
$$
C_{\ez}y^{1/\alpha-\ez}\leq b(y)\leq C_{\ez}'y^{1/\alpha+\ez}.
$$
By  (\ref{boundtransition}), we see
 \begin{eqnarray*}
 NP(B^{N,0}_{\ez_N^*}=0)=NP(B^{0}_{N\ez_N^*}=0)\leq C N
 b(N\ez_N^*)^{-d}\leq \frac{C}{C_{\ez}'}\frac{N (N\ez_N^*)^{d\ez}}
 {(N\ez_N^*)^{d/\alpha}}.
 \end{eqnarray*}
A simple calculation shows that given $\ez<1/2$, we can set
 \bgeqn
 \label{ezNstar}
\ez_N^*=N^{-\ez^*}\quad \textrm{ for }
\ez^*<1-\frac{\alpha}{d-\alpha d\ez}<1.
 \edeqn
Then $ NP(B^{N,0}_{\ez_N^*}=0)\rar0$ as $N\rar\infty$.
 When $d=\alpha=1$, since $s(x)\rar\infty$ as
$x\rar\infty$, we can set $x(0)=0$ and $\forall\, k\geq1$, there
exists $x(k)>x(k-1)$, such that if $x>x(k)$, $s(x)>k$. Then
$x(k)\rar\infty$ as $k\rar\infty$. Define function $s'$ on $\mbb
R^+$ such that $s'(x)=1$ for $0\leq x\leq x(1)$ and
$$s'(x)=k, \textrm{ for }x(k)< x\leq x(k+1) \text{ and } k\geq 1. $$
It is easy to see that $s'(x)\uparrow\infty$ as $x\rar\infty$ and
$\forall x>x(1)$, $s'(x)<s(x)$. Define
$$
\ez_N^*:=\left((\log N)\wedge\sqrt{s'(N/\log N)}\right)^{-1}.
$$
Then $N\ez_N^*\geq N/\log N$ and $N\ez_N^*\rar\infty$ as
$N\rar\infty$. Thus when $N$ is large enough ($N\ez_N^*>x(1)$),
 $$
 \ez_N^*s(N\ez_N^*)\geq s'(N\ez_N^*)/\sqrt{s'(N/\log N)}\geq \sqrt{s'(N/\log
 N)}\xrightarrow{N\rar\infty}\infty.
 $$
 We have
that
$$
NP(B^{N,0}_{\ez_N^*}=0)\leq C N
 b(N\ez_N^*)^{-1}=\frac{1}{\ez_N^*s(N \ez_N^*)}\rar0
$$
as $N\rar\infty$. Next,
 \beqnn
 \sum_{e\in\mbb S_N}p_N(e)P({\hat{\tau}}^N(\{0,e\})>\ez_N^*)
 &=&\sum_{e\in Z^d}p(e)P(\hat{\tau}(0,e)>N\ez_N^*)\\
 &\rar&\sum_{e\in Z^d}p(e)P(\hat{\tau}(0,e)=\infty)=\gamma_e.
 \eeqnn
Note that
$$P\left(\hat{\tau}^N(\{0,e\})\in(\ez_N^*,t]\right)
=P(\hat{\tau}^N(\{0,e\})>\ez_N^*)-P(\hat{\tau}^N(\{0,e\})>t).$$ Then
the second limit also holds. For any finite set $A\subset\mbb Z^d$,
$$
\sigma_N(A)=P(\hat{\tau}^N(A)\leq\ez_N^*)=P(\hat{\tau}(A)\leq
N\ez_N^*)\rar P(\hat{\tau}(A)<\infty)=\sigma(A).
$$
We are done. \qed

\smallskip

Next, we consider the `perturbation' term. As in \cite{[CP05]}, let
$P_F$ denote the set of finite subsets of $\mbb Z^d$. For $A\in
P_F$, $x\in\mbb S_N$, $\xi\in\{0,1\}^{\mbb S_N}$, define
$$
\chi_N(A,x,\xi)=\prod_{e\in A/{b(N)}}\xi(x+e).
$$
We also define
 \beqnn
 \beta_N(A)=\left\{\begin{array}{lll}
 \theta_0^N(p(e))^2,& A=\{e\},\\
 2\theta_0^Np(e)p(e'),&A=\{e,e'\},\\
 0,&\textrm{otherwise,}
 \end{array}\right.
 \eeqnn
and
 \beqnn
 \dz_N(A)=\left\{\begin{array}{lll}
 \theta_1^N,& A=\emptyset,\\
 \theta_1^N[(p(e))^2-2p(e)],&A=\{e\},\\
 2\theta_1^Np(e)p(e'),& A=\{e,e'\},\\
 0,&\textrm{otherwise}.
 \end{array}\right.
 \eeqnn

\begin{remark}
\label{Perburtion} According to the arguments in Section 1.2 of
\cite{[CP05]}, the `Perturbation assumptions' (P1) to (P5) there are
satisfied by the above coefficients with $l_N=b(N)$.
\end{remark}

The following proposition is exactly the same with Proposition 3.3
of \cite{[CP05]}. The Proposition 3.3 of \cite{[CP05]} was proved in
Section 4 there in which the proof of the results did not use any of
the kernel assumptions. Thus we can state the following proposition
without proof.

\smallskip

 \begin{proposition}
 \label{prop3.3}
 For $K,T>0$, there exists a finite constant $C_1(K,T)$ such that if
 $\sup_NX_0^N(1)\leq K$, then
 $$
 \sup_N E\left(\sup_{t\leq T}X_t^N(1)^2\right)\leq C_1(K,T).
 $$
 \end{proposition}
 \smallskip

This bound allows us to employ the $L^2$ arguments of \cite{[CP05]}.
Next, we consider another technical result, a version of Proposition
3.4 of \cite{[CP05]}. For $A\in P_F$, $\phi:[0,T]\times S_N\lra\mbb
R$ bounded and measurable, $K>0$ and $t\in[0,T]$, define
 \beqnn
 &&\mathcal{E}_N(A,\phi,K,t)\\
 &&\quad=\sup_{X_0^N(1)\leq K}E\left(\left(\int_0^t\left[
 \frac{1}{N}\sum_x\phi_s(x)\chi_N(A,x,\xi_s^N)-\sigma_N(A)X_s^N(\phi_s)
 \right]ds\right)^2\right).
 \eeqnn
Set $c_{\beta}=\sup_N|\theta_0^N|\sum_{e,e'\in\mbb Z^d}p(e)p(e')$
and $\bar{c}=c_{\beta}+k_{\dz}$, where $k_{\dz}=\sup_N|\theta_1^N|$.
The following proposition is a version of Proposition 3.4 of
\cite{[CP05]}.
 \begin{proposition}
 \label{Prop3.4}
 There is a positive sequence $\ez_N\lra0$ as $N\lra\infty$, and for
 any $K,T>0$, a constant $C_2(K,T)>0$, $\underline{\alpha}\leq1\wedge\alpha$,
 such that for any
 $\phi\in C_b([0,T]\times\mbb S_N)$ satisfying
 $\sup_{s\leq T}||\phi_s||_{\textrm{Lip}}\leq K$, nonempty $A\in
 P_F$, $\bar{a}\in A$, $J\geq1$ and $0\leq t\leq T$,
  \beqnn
 \mathcal{E}_N(A,\phi,K,t)\leq C_4(K,T)&\bigg{[}&\ez_N^*
 e^{\bar{c}\ez_N^*}+J^{-2}\\
&&\quad+J^2\left(\ez_N|A|+(\sigma_N(A)\wedge(\ez_N
+\left|\frac{\bar{a}}{b(N)}\right|^{\underline{\alpha}}))
\right)\bigg{]}. \eeqnn In particular,
$\lim_{N\rar\infty}\sup_{t\leq T}\mathcal{E}_N(A,\phi,K,t)=0.$
 \end{proposition}
 {\bf Proof. } We can follow the arguments in Section 5 and Section
 6 of \cite{[CP05]}. In fact, only a small trick is needed.
 For $\alpha\in(0,2]$ and $d>\alpha$, we may find an $\underline{\alpha}<\alpha$
 which is close enough to $\alpha$ so that
 \bgeqn
 \label{Nezlimit}
 E(|B_{\ez_N^*}^{N,0}|^{\underline{\alpha}})
 =\frac{N\ez_N^*|p|_{\underline{\alpha}}}{b(N)^{\underline{\alpha}}}
 \lra0\quad\textrm{ as } N\lra\infty.
  \edeqn
(Note that $b$ is a function of regular variation with index
$1/\alpha$ and recall the choice of $\ez_N^*$ in Lemma
\ref{Kerlemma} when $d>\alpha$). Fix this $\underline{\alpha}$.
 For $||\phi||_{\textrm{Lip}}\leq K$, (\ref{LIP}) implies
 \beqnn
 E\left(\left|\phi\left(y-\frac{\bar{a}}{b(N)}+B_s^{N,0}\right)
 -\phi(y)\right|\right)&\leq& 2K
 E\left|B_s^{N,0}-\frac{\bar{a}}{b(N)}\right|^{\underline{\alpha}\wedge1}\\
 &\leq&2KE\left(|B_{s}^{N,0}|^{\underline{\alpha}\wedge1}\right)
  +2K\left|\frac{\bar{a}}{b(N)}\right|^{\underline{\alpha}\wedge1}.\eeqnn
When $\alpha>1$, we may assume
$\underline{\alpha}\wedge1=1\leq\underline{\alpha}$.
(\ref{Nezlimit}) suggests
$$
E\left(|B_{\ez_N^*}^{N,0}|^{\underline{\alpha}\wedge1}\right)\lra0\quad\textrm{
as } N\lra\infty.
$$
When $d=\alpha=1$, for any $\underline{\alpha}<1$, by (\ref{LIP}),
\beqnn
 && E\left(\left|\phi\left(y-\frac{\bar{a}}{b(N)}+B_s^{N,0}\right)
 -\phi(y)\right|\right)\\
 &&\quad\leq 2K
 E\left(|B_s^{N,0}-{\bar{a}}/{b(N)}|^{\underline{\alpha}}; |B_s^{N,0}|<s(N)^{-1}
 \right)
 +2||\phi||_{\infty}P\left(|B_s^{N,0}|\geq s(N)^{-1}\right)\\
 &&\quad\leq \frac{2K}{s(N)^{\underline{\alpha}}}
 +2K\left|\frac{\bar{a}}{b(N)}\right|^{\underline{\alpha}}
   +2KP\left(|B_s^{N,0}|\geq s(N)^{-1}\right)
 .\eeqnn
We want to estimate the last term above for $s=\ez_N^*$. First,
$$P\left(|B_{\ez_N^*}^{N,0}|\geq s(N)^{-1}\right)
=P\left(|B_{N\ez_N^*}^{0}|\geq N\right).$$
  By Proposition \ref{propincrease} and (\ref{regular}),
$P(|B_{\ez_N^*}^{N,0}|\geq s(N)^{-1})$ is bounded by
$$
C_{\ref{increasein}}\frac{N\ez_N^*}{l(N)}=
C_{\ref{increasein}}\frac{l(N\ez_N^*s(N\ez_N^*))}{l(N)}\leq
\frac{C_{\ref{increasein}}s}{C_{\ez}(\ez_N^*s(N\ez_N^*))^{1-\ez}}.
$$
Recall the choice of $\ez_N^*$ in the Lemma \ref{Kerlemma} when
$d=\alpha=1$. The last term above goes to zero when $N\rar\infty$.
Set
$$\ez_N=2KE(|B_{\ez_N^*}^{N,0}|^{\underline{\alpha}\wedge1})\quad\textrm{ for
}d>\alpha$$ and
$$\ez_N=\frac{2K}{s(N)^{\underline{\alpha}}}
   +2KP(|B_{\ez_N^*}^{N,0}|\geq s(N)^{-1})\quad\textrm{ for }d=\alpha=1.$$
Then $\ez_N\rar0$ as $N\rar\infty$ and
\beqlb\label{keyineq}E\left(\left|\phi\left(y-\frac{\bar{a}}{b(N)}
+B_{\ez_N^*}^{N,0}\right)
 -\phi(y)\right|\right)\leq \ez_N
  +2K\left|\frac{\bar{a}}{b(N)}\right|^{\underline{\alpha}\wedge1}.\eeqlb
With (\ref{keyineq}) in mind, the reader may go back to
\cite{[CP05]} for the proof of this proposition. In fact, as in
\cite{[CP05]}, we first define $\eta_N$ as (5.1) of \cite{[CP05]}
and  decompose it into four error terms $\eta_i^N, i=1,2,3,4$. And
decompose $\eta_3^N$ into two terms, $\eta_{3,1}^N$ and
$\eta_{3,2}^N$, as in (5.15) and (5.16) of \cite{[CP05]}
respectively. (\ref{keyineq}) will be used when we estimate
$\eta_{3,2}^N(s)$ as on p.944 of \cite{[CP05]}. Only a part of the
proof at the end of Section 5 of \cite{[CP05]} is needed to be
modified. When estimate $\eta_{3,1}^N$, we also need (\ref{A3}).\qed

The following technical lemma will be used in checking the Compact
Containment Condition.

\begin{lemma}
\label{CCCtech}Let $P_t^N$ denote the semigroup associated with
generator ${\cal A}_N$. We have
$$
X_0^N\left(P_s^N(1_{B(0,n)^c})\right)\lra0\quad \textrm{ as }
n\lra\infty
$$ uniformly in  $N$  and $s\leq t$.
\end{lemma}
{\bf Proof. }  Since
$$
X_0^N\left(P_s^N(1_{B(0,n)^c})\right)\leq
X_0^N\left(B(0,n/2)^c\right)+X_0^N(1)P(|B_s^{N,0}|>n/2),
$$
and (A2) holds, it suffices to show $P(|B_s^{N,0}|>n/2)$ goes to 0
uniformly as $n\rar\infty$.  For $0<c<1$, note that
 \bgeqn
 \label{for1}P(|B_s^{N,0}|>cn)=P(|B_{Ns}^0|>cnb(N)).
 \edeqn
 When $\alpha=2$, the desired result follows from Chebyshev's
 inequality. We only need to consider the case of $\alpha<2$.
 Clearly, we can deal
separately with the different coordinates of $B_s^{N,0}$ and the
distribution of each coordinate of $Y_1$ is  a dimension-one
$(\sigma^2,\alpha)$-stable distribution. (A1) implies that each
coordinate of $p(\cdot)$ is in the domain of attraction of the
dimension-one $(\sigma^2,\alpha)$-stable distribution. Thus, for
this proof only,
 we can assume $d=1$ (Here we drop the assumption $d\geq \alpha$).
 By Proposition \ref{propincrease} and
(\ref{regular}), the right hand side of (\ref{for1}) is bounded by
$$
C_{\ref{increasein}}\frac{Ns}{l(cnb(N))}=
C_{\ref{increasein}}\frac{l(b(N))s}{l(cnb(N))} \leq
\frac{C_{\ref{increasein}}s}{C_{\ez}(cn)^{\alpha-\ez}},
$$
where the inequality holds for $cn>1$. The desired result is then
immediate.  \qed

\noindent {\bf Proof of Theorem \ref{mainUP}.} Now, we are in
position to prove Theorem \ref{mainUP}. First, we check the compact
containment condition. Let $h_n:\mbb R^d\rar[0,1]$ be a $C^{\infty}$
function such that
 $$
 B(0,n)\subset\{x:h_n(x)=0\}\subset\{x:h_n(x)<1\}\subset B(0,n+1)
 $$
and
 $$
 \sup_n\sum_{i,j,k\leq d}||(h_n)_i||_{\infty}+||(h_n)_{ij}||_{\infty}
 +||(h_n)_{ijk}||_{\infty}\equiv C_h<\infty. $$
 Let $\phi_n=\sigma^2|\Delta^{\alpha/2}h_n|/2$.
Using Taylor's formula and dominated convergence theorem we obtain
there exists a constant $C>0$ such that
$$\sup_n\sum_{i\leq d}||(\Delta^{\alpha/2}h_n)_i||_{\infty}<C.$$
Thus $\sup_n||\phi_n||_{\textrm{Lip}}<C$. We may define $\dz_N^1$,
$\dz_N^2$ and $d_0^N$ as on p.927 of \cite{[CP05]}. With Proposition
\ref{Prop3.4} in hand one can check that  both Lemma 3.6 and
Proposition 3.8 in \cite{[CP05]} are available. To establish the
Compact Containment Condition, we may follow the proof of
Proposition 3.9 of  \cite{[CP05]}. In fact, the argument above and
Lemma \ref{CCCtech} show that
$$
\lim_{(N,n)\rar\infty}E\left( \int_0^tX_s^N (|{\cal A}_Nh_n|
)ds\right)=0.
$$
Then the following argument for the compact containment condition
are exactly the same with that in \cite{[CP05]}. Next, with Lemma
\ref{conge}, Proposition \ref{prop3.3} and Lemma \ref{estimateMNt}
in hand, the proof of C-tightness is analogous to that of
Proposition 3.7 of \cite{[CP05]}. By Proposition \ref{prop3.3}, we
see that the $L^2$-method in \cite{[CP05]} is available. Thus, we
may use the arguments in the proof of Proposition 3.2 in
\cite{[CP05]} with some trivial modifications to obtain the desired
convergence theorem, Theorem \ref{mainUP}. \qed

\section{Proof of Theorem \ref{main2}}
In this section we assume that
$$
d=\alpha=1\qquad\textrm{ and }\qquad b(t)=t.
$$
With $\underline{\alpha}$ we always mean a constant which is
strictly less than 1. We can adopt some of the arguments of
\cite{[CP08]} to prove some analogous results to those in
\cite{[CP08]} without using the fact that $p(\cdot)$ is in the
domain of attraction of a stable law. We will refer the reader to
these results as we use them.

\subsection{Characterization of $\gamma^*$}
Recall the definitions of $\hat{\tau}$ and $\tau$ in Section
\ref{MR}. For $e,e'\in \mbb Z$ define the event
$\Gamma_T(e,e')=\{\hat{\tau}(e,e')<T,\hat{\tau}(0,e)\wedge
\hat{\tau}(0,e')>T\}$, and let

\bgeqn \label{QT} q_T=\sum_{e,e'}p(e)p(e')P(\Gamma_T(e,e')).
 \edeqn
We have the following characterization of $\gamma^*$.
\begin{proposition}
\label{Prop2.1} \bgeqn\label{gamma*} \gamma*=\lim_{T\rar\infty}
(\log T)q_T<\infty. \edeqn
\end{proposition}
To prove Proposition \ref{Prop2.1}, we follow the arguments in
Section 2 of \cite{[CP08]}. Let $\tau_x=\inf\{t\geq0:B_t^0=x\}$, and
write $P^x$ to indicate the law of the walk $B^x_{\cdot}$. Let
$\tilde{P}(\cdot)=\sum_ep(e)P^e(\cdot)$, and define
 \bgeqn
 \label{Ht}
H(t)=\tilde{P}(\tau_0>t).
 \edeqn
The following proposition is a version of Proposition 2.2 of
\cite{[CP08]}.
 \begin{proposition}
 \label{prop2.2}
 \bgeqn
 \label{2.4}
 \lim_{t\rar\infty}H(t)\log t=p_1(0)^{-1}.
 \edeqn
 \bgeqn
 \label{2.5}
 \frac{P^x(\tau_0>t)}{H(t)}\leq 2a(x)\quad\textrm{for all}\quad
 x\in\mbb Z,\,t>0.
 \edeqn
 \bgeqn
 \label{2.6}
 \lim_{t\rar\infty}\frac{P^x(\tau_0>t)}{H(t)}=a(x)\quad\textrm{for
 all}\quad x\in \mbb Z.
 \edeqn
 \bgeqn
 \label{2.7}
 a(x)/|x|,\, x\neq0\textrm{ is bounded on }\mbb Z.
 \edeqn
 \end{proposition}
{\bf Proof. }For (\ref{2.4}), let $G(t)=\int_0^tp_s(0,0)ds$.
Proposition \ref{CPtostable} implies $G(t)\sim p_1(0)\log t$ as
$t\rar\infty$ in $d=1$. Then one can follow the arguments in the
proof of Lemma A.3 in \cite{[CDP01]} by using the last exit time
decomposition of Lemma A.2 there and with (A.7)  replaced by
(\ref{boundtransition}) to obtain that $G(t)H(t)\rar1$ as
$t\rar\infty$; see the arguments after (A.8) of \cite{[CDP01]}. Then
(\ref{2.4}) holds.

Recall that $\{Z_n:n=0,1,2,\cdots\}$ denote the discrete time stable
random walk defined in Section \ref{ourmodel}. With abuse of
notation, let $P^x$ denote the law of the walk starting at $Z_0=x$.
Let $\sigma_x=\inf\{n\geq1:Z_n=x\}$. By T29.1 of \cite{[Sp76]},
$$
a(x)=\lim_{n\rar\infty}\sum_{k=0}^n[P^0(Z_k=0)-P^0(Z_k=x)]<\infty
\quad\textrm{exists for all }x\textrm{ in }\mbb Z.
$$
Note that P11.1, P11.2 and P11.3 in Chapter III of \cite{[Sp76]} are
available for one-dimensional recurrent random walk; see arguments
before P28.1 of \cite{[Sp76]}. Meanwhile, according to T29.1 and
P30.1 of \cite{[Sp76]}, (i)' and (ii)' on page 116 in Chapter III of
\cite{[Sp76]} also hold for one-dimensional random walk. Then we can
check that both P11.4 and P11.5 in Chapter III of \cite{[Sp76]} are
also available. Thus we have $$P^0(\sigma_x<\sigma_0)=1/2a(x).$$
Since the sequences of states visited by the walk $B_t^0$ is equal
in law to the sequences visited by the walk $Y_n$ (with $Y_0=0$), we
have $\tilde{P}(\tau_x<\tau_0)=1/2a(x)$. The strong Markov property
implies that
$$
H(t)\geq\sum_ep(e)P^e(\tau_x<\tau_0,\tau_0>t)\geq\sum_e
P^e(\tau_x<\tau_0)P^x(\tau_0>t)
$$
and then (\ref{2.5}) follows.

For (\ref{2.6}), by T32.1 of \cite{[Sp76]},
 \bgeqn
 \label{2.8}
 \lim_{n\rar\infty}\frac{P^x(\sigma_0>n)}{P^0(\sigma_0>n)}=a(x).
 \edeqn
Define
$$
h(n)=\sum_{0\leq k\leq n}P^0(Y_k=0).
$$
Then
 \bgeqn
 \label{asmhn}
h(n)\sim p_1(0)\sum_{k=1}^n\frac{1}{k} \quad\textrm{as }n\rar\infty;
 \edeqn
see Page 696 of \cite{[LR91]}. We also have that
$$
P^0(\sigma_0>n)=\frac{1}{h(n)}+o\left(\frac{1}{h(n)^2}\right);$$
 see the proof of Theorem 6.9 of \cite{[LR91]}. Thus
 \bgeqn
 \label{2.9}
 P^0(\sigma_0>n)\log n\rar p_1(0)^{-1}.
 \edeqn
   According to  a standard large deviations estimate for a rate-1 Poisson
process, say $S(t)$,  $e^{Ct}P(S(t) \notin[t/2,2t])\rar0$ as
$n\rar\infty$ for a some constant $C>0$. Then the fact that
$Y_{S(\cdot)}$ is a realization of $B^0_{\cdot}$ yields
$$
(1-o(e^{-Ct}))P^x(\sigma_0>2t)\leq P^x(\tau_0>t)\leq
o(e^{-Ct})+P^x(\sigma_0>t/2).
$$
The inequalities above, together with (\ref{2.8}) and (\ref{2.9}),
imply
 \bgeqn
 \label{2.10}
\lim_{t\rar\infty}\frac{P^x(\tau_0>t)}{P^x(\sigma_0>t)}=1.
 \edeqn
By (\ref{2.4}) we see  $H(t)/P^0(\sigma_0>t)\rar1$
 as $t\rar\infty$. Then (\ref{2.8}) and (\ref{2.10})
 tell us (\ref{2.6}) holds readily.
Finally, (\ref{2.7}) follows from the fact that
$$\lim_{|x|\rar\infty}\frac{a(x)}{|x|}=0;$$ see P29.3 of
\cite{[Sp76]} and elsewhere. We have completed the proof. \qed

  The proof of Proposition \ref{Prop2.1} is now exactly
 as that of Proposition 2.1 in Section 2 of \cite{[CP08]}.
 We omit it here.

\subsection{Voter and Biased Voter  Estimates}
 In this subsection, we consider voter, biased voter
bounds. We follow the arguments in Section 5 of \cite{[CP08]} step
by step. For $b,\nu\geq0$, the 1-biased voter model $\bar{\xi}_t$ is
the Feller process taking values in $\{0,1\}^{\mbb Z}$, with rate
function
 \beqlb
 \label{5.1}
 \bar{c}(x,\xi)=\begin{cases}
 (\nu+b)f_1(x,\xi)& \textrm{ if }\xi(x)=0,\cr
 \nu f_0(x,\xi)& \textrm{ if }\xi(x)=1,
                      \end{cases}
 \eeqlb
where $f_i(x,\xi)$ is as in (\ref{1.1}). The 0-biased voter model is
the Feller process $\underline{\xi}_t$  taking values in
$\{0,1\}^{\mbb Z}$ with rate function
  \beqlb
 \label{5.2}
 \underline{c}(x,\xi)=\begin{cases}
 \nu f_1(x,\xi)& \textrm{ if }\xi(x)=0,\cr
 (\nu+b) f_0(x,\xi)& \textrm{ if }\xi(x)=1.
                      \end{cases}
 \eeqlb
The voter model $\hat{\xi}_t$ is the 1-biased voter model with bias
$b=0$. Then by Theorem III.1.5 of \cite{[L85]}, assuming
$\underline{\xi}_0=\hat{\xi}_0=\bar{\xi}_0$, we may define
$\underline{\xi}_t, \hat{\xi}_t$ and $\bar{\xi}_t$ on a common
probability space so that \bgeqn
\label{5.3}\underline{\xi}_t\leq\hat{\xi}_t\leq\bar{\xi}_t\textrm{
for all }t\geq0.\edeqn For $\xi,\zeta\in\{0,1\}^{\mbb Z}$,
$\xi\leq\zeta$ means $\xi(x)\leq\zeta(x)$ for all $x\in \mbb Z$.

Let us recall the voter model duality; see \cite{[L85]}. Recall also
the coalescing random walk system $\{\hat{B}_t^x:x\in \mbb Z\}$
defined in Subsection \ref{MR}. The duality equation for the rate-1
($\nu=1$) voter model is: for finite $A\subset\mbb Z$,
 \bgeqn
 \label{5.4}
 P(\hat{\xi}_t(x)=1\forall x\in A)=P(\hat{\xi}_0(\hat{B}_t^x)=1\forall x\in
 A).
 \edeqn
Define the mean range of the random walk $B_t^0$ by
$$
R(t)=E\left(\sum_x1_{\{B_s^0=x \textrm{ for some } s\leq
t\}}\right).$$ By a result for the range of the discrete time stable
random walk in \cite{[LR91]},
 \bgeqn
 \label{5.5}
\lim_{t\rar\infty}\frac{R(t)}{t/\log t}=p_1(0)^{-1};
 \edeqn
see (1.e) of \cite{[LR91]} and recall  (\ref{asmhn}) for the
asymptotic behavior of $h(n)$.

First, we consider the voter estimates. Let $P_t,t\geq0$ be the
semigroup of a rate-1 random walk with step distribution $p(\cdot)$.
Recall the definition of $|p|_{\underline{\alpha}}$ in Section 3.
For $\phi:\mbb Z\rar\mbb R$ and $\xi\in\{0,1\}^{\mbb Z}$, let
$$\xi(\phi)=\sum_x\phi(x)\xi(x).$$
\begin{lemma}
\label{lemma5.1} Let $\hat{\xi}_t$ denote the rate-$\nu$ voter
model. Then for all bounded $\phi:\mbb Z\rar\mbb R^+$,
$0<\underline{\alpha}<1$ and $t\geq0$,
 \bgeqn\label{5.8}
 E(\hat{\xi}(\phi f_0(\hat{\xi}_t)))\leq(\nu t|p|_{\underline{\alpha}}H(2\nu
 t))^{1/2}||\phi||_{\underline{\alpha}/2}|\bar{\xi}_0|+H(2\nu
 t)\hat{\xi}_0(\phi).
  \edeqn
 \end{lemma}
\begin{remark}
\label{remark4.1}
 (\ref{5.8}) is just a version of (5.8) in Lemma 5.1 of
 \cite{[CP08]}. We slightly abuse our notation and we can prove
  that the other statements in Lemma 5.1  of \cite{[CP08]}
  ((5.6), (5.7) and (5.9) there) hold without modifying any
  arguments of their proofs.
\end{remark}
\begin{remark}
Recall the definition of $||\phi||_{\underline{\alpha}}$ in Section
3. We see for $\phi=1$, the right side of (\ref{5.8}) is just
$H(2\nu t)|\hat{\xi}_0|$.
\end{remark}
 {\bf Proof. } It suffices to consider $\nu=1$.
 Using the voter duality equation (\ref{5.4})
 and following the arguments in the proof of (5.8) of
 \cite{[CP08]}, we have
 $$E(\hat{\xi}(\phi
 f_0(\hat{\xi}_t)))\leq\sum_{e,z}\hat{\xi}_0(z)p(e)
 E\left(\phi(z+B_t^0)1_{\{\tau(0,e)>t\}}\right).
 $$
For any $z$ and $0<\underline{\alpha}<1$,
 \beqnn
 &&\sum_ep(e)E\left(\phi(z+B_t^0)1_{\{\tau(0,e)>t\}}\right)\cr
 &&\quad\quad\leq\sum_ep(e)
 E\left(\left(||\phi||_{\underline{\alpha}/2}|B_t^0|^{\underline{\alpha}/2}+\phi(z)\right)
 1_{\{\tau(0,e)>t\}}\right)\cr
 &&\quad\quad\leq
 ||\phi||_{\underline{\alpha}/2}\left(E(|B_t^0|^{\underline{\alpha}})
 \sum_ep(e)P(\tau(0,e)>t)\right)^{1/2}\cr
 &&\quad\quad\quad+\phi(z)\sum_ep(e)P(\tau(0,e)>t).
 \eeqnn
Since $E(|B_t^0|^{\underline{\alpha}})\leq
t|p|_{\underline{\alpha}},$ this proves (\ref{5.8}).  \qed

Next, we give some biased voter model bounds. Let $\bar{\xi}_t$ be
the 1-biased voter model with rate function (\ref{5.1}). By the same
arguments in Section 4 of \cite{[CP05]}, we can prove the following
inequalities without using any of kernel assumptions.
 \beqlb
 \label{5.12}
 E(|\bar{\xi}_t|)&\leq& e^{bt}|\bar{\xi}_0|,\\
 \label{5.13}
 E(|\bar{\xi}_t|^2)&\leq& e^{2bt}\left(|\bar{\xi}_0|^2
 +\frac{2\nu+b}{b}(1-e^{-bt})|\bar{\xi}_0|\right)\\
 \label{5.14}
 &\leq& e^{2bt}\left(|\bar{\xi}_0|^2
 +(2\nu+b)t|\bar{\xi}_0|\right)
 \eeqlb
In the subsection 4.3 below, we will compare the Lotka-Volterra
model $\xi_t^N$ with the biased voter models
$\underline{\xi}_t^N,\bar{\xi}_t^N$ on $\mbb S_N$. In order to
construct coupling $\underline{\xi}_t^N\leq \xi_t^N\leq
\bar{\xi}_t^N$ we assume that the voting and bias rates $\nu_N$ and
$b_N$ are
 \bgeqn
 \label{5.15}
 \nu=\nu_N=N-\bar{\theta}\log N \textrm{ and
 }b=b_N=2\bar{\theta}\log N.
 \edeqn
As in \cite{[CP08]}, we need  improved versions of (\ref{5.12}) and
(\ref{5.13}). For $p\geq2$ and $0<\underline{\alpha}<1$ define
 \beqnn
 \kappa_p&=&\kappa_p(b,\nu)=3(b H(2\nu/b^p)+e^2)\textrm{ and
 }\kappa=\kappa_3,\cr
 A&=& A(b,\nu)=bR(2\nu/b^3)+3e^2(1+2\nu/b),\cr
 B_p&=& B_p(b,\nu,\underline{\alpha})=(|p|_{\underline{\alpha}}\nu
 b^{2-p}H(2\nu/b^p))^{1/2}+b
 H(2\nu/b^p)(|p|_{\underline{\alpha}}(\nu/b^p+1))^{1/2}
 \eeqnn
and
 \beqnn
 h_1(b,\nu)(t)&=& e^2t^{-1/3}+2\kappa e^{2+2\kappa t},\cr
 h_2(b,\nu)(t)&=& e^2t^{-1/3}(1+2\nu/b)+5\kappa A e^{1+3\kappa
 t}.
 \eeqnn
Put $P\phi(x)=\sum_yp(y-x)\phi(y)$ and define the operators \bgeqn
\label{defA*} \bar{\cal A}\phi=\nu(P\phi-\phi)\textrm{ and }{\cal
A}^*=(1+b/\nu)\bar{\cal A} \edeqn and denote the associated
semigroups by $\bar{P}_t$ and  $P^*_t$ respectively.
 \begin{remark}
 Comparing the constants and functions defined above with those
 defined in (5.16) and (5.17) of \cite{[CP08]}, we see that only $B_p$ is
 different. We replaced $2\sigma^2$ by $|p|_{\underline{\alpha}}$.
 \end{remark}
 \begin{remark}
 \label{remark5.3}
 For the parameters $\nu=\nu_N$, $b=b_N$ in (\ref{5.15}),
 (\ref{2.4}) and (\ref{5.5}) imply that $\kappa_p=O(1)$, $A=O(N/\log
 N)$ and $B_p=O(N^{1/2}(\log N)^{(1-p)/2})$ as $N\rar\infty$.
 \end{remark}
  \begin{remark}
 \label{LVreasoning}
 The estimates in Remark \ref{remark5.3} will play important roles
 in the following proofs. That is why we are forced to assume that
 $\{p(x)\}$ is in the  domain of normal attraction of a stable law. Or
 we need to replace $\log N$ by $\int_1^N b(s)^{-1}ds$. Then the
 estimates in Remark \ref{remark5.3} will be not available.
 \end{remark}
The following proposition is a version of Proposition 5.4 of
\cite{[CP08]}.
 \begin{proposition}
 \label{prop5.4}
 Assume $b\geq1$ and $p\geq2$. For all $t\geq0$,
 \beqlb
 \label{5.19}
 E(|\bar{\xi}_t|)&\leq& e^{b^{1-p}+\kappa_pt}|\bar{\xi}_0|,\\
 \label{5.20}
 E(|\bar{\xi}_t|^2)&\leq& e^{2+2\kappa t}|\bar{\xi}_0|^2+4Ae^{1+3\kappa
 t}|\bar{\xi}_0|,\\
 \label{5.21}
 b E(\bar{\xi}_t(f_0(\bar{\xi}_t)))&\leq& h_1(t)|\bar{\xi}_0|,\\
 \label{5.22}
 b E(|\bar{\xi}_t|\bar{\xi}_t(f_0(\bar{\xi}_t)))
 &\leq& h_1(t)|\bar{\xi}_0|^2+h_2(t)|\bar{\xi}_0|.
 \eeqlb
 For all bounded $\phi:\mbb Z\rar[0,\infty)$, $p\geq3$ and
 $0<\underline{\alpha}<1$,
 \beqlb
 \label{5.23}
 E(\bar{\xi}_t(\phi))\leq e^{b^{1-p}+(1+\kappa_p)t}
 \left(\bar{\xi}_0(P_t^*(\phi))+\left[\kappa_pb^{2-p}
 ||\phi||_{\infty}+B_p||\phi||_{\underline{\alpha}/2}
 \right]|\bar{\xi}_0| \right).
 \eeqlb
 \end{proposition}
\begin{remark}
\label{remark4.5} Proposition 5.4 of \cite{[CP08]} was proved with
the help of Lemma 5.1, Lemma 5.5 and Lemma 5.6 there. We can adopt
the arguments in \cite{[CP08]} to obtain similar results in Lemma
5.5 and Lemma 5.6 of \cite{[CP08]}. With abuse of notation, in the
following we assume that
 those two lemmas are available for us.
\end{remark}
\begin{remark}
The only difference between Proposition \ref{prop5.4} and
Proposition 5.4 of \cite{[CP08]} is that inequality (\ref{5.23}) is
different from inequality (5.23) there. In fact, the key reason is
that when prove the inequality (\ref{5.23}), we will use estimate
(\ref{5.8}) in Lemma \ref{lemma5.1} of this paper replacing the
estimate (5.8) of Lemma 5.1 of \cite{[CP08]}.
\end{remark}
{\bf Proof. }According to Remark \ref{remark4.1}, Remark
\ref{remark4.5} and the coupling (\ref{5.3}),  we can follow the
arguments in \cite{[CP08]} to obtain that (5.36), (5.37) and (5.38)
there are available which will be used in the following proof.  Put
$\ez=b^{-p}$ and assume $\phi\geq0$. We also have that
 \beqlb
 \label{5.39}
 E(|\bar{\xi}_{\ez}(b\phi f_0(\bar{\xi}_{\ez}))&-&\hat{\xi}_{\ez}(b\phi
 f_0(\hat{\xi}_{\ez}))|)\cr
 &&\leq2b||\phi||_{\infty}E(|\bar{\xi}_{\ez}|-|\hat{\xi}_{\ez}|)\leq
 2b(e^{b\ez}-1)||\phi||_{\infty}|\bar{\xi}_0|
 \eeqlb
which is just a version of (5.39) of \cite{[CP08]} (In fact, they
are the same). The voter model estimate (\ref{5.8}) tells us
 \beqlb
 \label{5.40}
 E(\bar{\xi}_{\ez}(b\phi f_0(\bar{\xi}_{\ez})))&\leq&
 2eb^2\ez||\phi||_{\infty}|\bar{\xi}_0|\cr
 &&+b(|p|_{\underline{\alpha}}\nu\ez
 H(2\nu\ez))^{1/2}||\phi||_{\underline{\alpha}/2}|\bar{\xi}_0|+b
 H(2\nu\ez)\bar{\xi}_0(\phi).
 \eeqlb
 By using Markov property, we see for $s\geq\ez$,
 \beqlb
 \label{5.41}
 &&
 E(\bar{\xi}_{s}(b\phi f_0(\bar{\xi}_{s}))|{\cal
 F}_{s-\ez})\cr
  &&\quad
  \leq\left(2eb^2\ez||\phi||_{\infty}+b(|p|_{\underline{\alpha}}\nu\ez
 H(2\nu\ez))^{1/2}||\phi||_{\underline{\alpha}/2}\right)|\bar{\xi}_{s-\ez}|+b
 H(2\nu\ez)\bar{\xi}_{s-\ez}(\phi).
 \eeqlb
Take expectations in (\ref{5.41}) for $\phi=1$ and recall the
definition $||\phi||_{\underline{\alpha}}$ in Section 3. We have for
$s\geq\ez$
 \beqlb
 \label{5.42}
 E(\bar{\xi}_{s}(b\phi f_0(\bar{\xi}_{s})))\leq
 \kappa_pE(|\bar{\xi}_{s-\ez}|).
 \eeqlb
Using this inequality in (5.36) of \cite{[CP08]} yields for
$s\geq\ez$,
$$
E(|\bar{\xi}_t|)\leq E(|\bar{\xi}_{\ez}|)
+\kappa_p\int_{\ez}^tE(|\bar{\xi}_{s-\ez}|)ds\leq
e^{b\ez}+\kappa_p\int_0^tE(|\bar{\xi}_s|)ds,
$$
where the second inequality follows from (5.38) of \cite{[CP08]}.
This bound also holds for $t\leq \ez$. Then Gronwall's inequality
implies that (\ref{5.19}) holds.

Again using (5.38) of \cite{[CP08]} gives that for $\psi:\mbb
Z\rar\mbb R^+$,
 $$
 |E(\bar{\xi}_{\ez}(\psi))-\bar{\xi}_0(\psi)|\leq
 (e^{b\ez}-1)\bar{\xi}_0(P^*_{\ez}\psi)+|\bar{\xi}_0(P^*_{\ez})
 -\bar{\xi}_0(\psi)|.
 $$
Note that
$$
 |P_{\ez}^*\psi(x)-\psi(x)|\leq||\psi||_{\underline{\alpha}/2}
 E(|B^0_{\nu\ez(1+b/\nu)}|^{\underline{\alpha}/2})
 \leq||\psi||_{\underline{\alpha}/2}(\ez(\nu+b)|p|_{\underline{\alpha}})^{1/2}.
$$
Thus $$|E(\bar{\xi}_{\ez}(\psi))-\bar{\xi}_0(\psi)|\leq
\left(eb\ez||\psi||_{\infty}+
||\psi||_{\underline{\alpha}/2}(\ez(\nu+b)|p|_{\underline{\alpha}})^{1/2}\right)|\bar{\xi}_0|.
$$
Then by using Markov property, for $s\geq\ez$,
$$
E(\bar{\xi}_{s-\ez}(\psi))\leq E(\bar{\xi}_s(\psi))+
\left(eb\ez||\psi||_{\infty}+
||\psi||_{\underline{\alpha}/2}(\ez(\nu+b)|p|_{\underline{\alpha}})^{1/2}\right)
E(|\bar{\xi}_{s-\ez}|).
$$
Since
$||P_{t-s}^*\phi||_{\underline{\alpha}/2}\leq||\phi||_{\underline{\alpha}/2}$,
using above inequality in (\ref{5.41}) with $\psi=P_{t-s}^*\phi$
replacing $\phi$, we have for $s\geq\ez$,
 \beqlb
 \label{5.43}
 E(\bar{\xi}_{s}(bP^*_{t-s}\phi f_0(\bar{\xi}_{s})))
  \leq\left(\kappa_p b^2\ez||\phi||_{\infty}+B_p||\phi||_{\underline{\alpha}/2}\right)
  E(|\bar{\xi}_{s-\ez}|)+\kappa_p E(\bar{\xi}_{s}(P^*_{t-s}\phi)),
 \eeqlb
 which is a version of (5.43) of \cite{[CP08]}.
 Then the following arguments for proving (\ref{5.23}) are very similar to
 those after (5.43) in
 \cite{[CP08]}.
 We have proved (\ref{5.19}) and (\ref{5.23}). The other
 statements in the proposition can be proved in a similar way to that used to prove
 their counterparts in \cite{[CP08]} (recall  Remark \ref{remark4.1},
  Remark \ref{remark4.5}). We omit it here. \qed
\begin{remark}
\label{remarksec5} We have followed the arguments in Section 5 of
\cite{[CP08]} to obtain some voter and biased voter estimates. In
fact, we only replaced (5.8) and (5.23) in Section 5 of
\cite{[CP08]} by (\ref{5.8}) and (\ref{5.23}) respectively and
modified the arguments in the proof of (5.19) and (5.23) of
\cite{[CP08]}; please compare (\ref{5.40})-(\ref{5.43}) with their
counterparts (5.40)-(5.43) in Section 5 of \cite{[CP08]}. We can
also adopt the arguments there to obtain similar results to  all
other statements in Section 5 of \cite{[CP08]} without using the
fact the $p(\cdot)$ is in the domain of attraction of a stable law.
In the next subsection, we will directly refer to them.
\end{remark}

\subsection{Four Key Results}
In this subsection, we will give analogous results to Propositions
4.3, 4.4, 4.5 and 4.7 of \cite{[CP08]}. We first list those results
and will give their proofs later. Let \bgeqn \label{4.1}
g(s)=C_{\ref{4.1}}s^{-1/3}e^{C_{\ref{4.1}}s},\edeqn where
$C_{\ref{4.1}}$ will be chosen later.
 \begin{proposition}
 \label{prop4.3}
 (a) For $T>0$ there is a constant $C_{\ref{4.2}}(T)$ such that for all $N\in \mbb
 N$,
 \beqlb
 \label{4.2}
 \sup_{t\leq T}E(X_t^N(1))&\leq& C_{\ref{4.2}}(T)X_0^N(1),\\
 \label{4.3}
 E\left(\sup_{t\leq T}X_t^N(1)^2\right)&\leq&
 C_{\ref{4.2}}(T)(X_0^N(1)^2+X_0^N(1)).
 \eeqlb
 (b) For all $s>0$ and $N\in\mbb N$,
 \beqlb
 \label{4.4}
 (\log N)E(X_s^N(f_0^N(\cdot,\xi_s^N)))&\leq& g(s)X_0^N(1),\\
 \label{4.5}
 (\log N)E(X_s^N(1)X_s^N(f_0^N(\cdot,\xi_s^N)))&\leq&
 g(s)(X_0^N(1)^2+X_s^N(1)).
 \eeqlb
 \end{proposition}
Let ${\cal A}^*_N(\psi)=\frac{1}{N}(N+\bar{\theta}\log N){\cal
A}_N(\psi)$ with semigroup $P_t^{N,*}$.
 \begin{proposition}
 \label{prop4.4}
 For $p\geq 3$ there is a constant $C_{\ref{4.6}}(p)$ such that for
 any $t\geq0$ and $\phi:\mbb R\rar\mbb R^+$,
 \beqlb
 \label{4.6}
 E(X_t^N(\phi))&\leq& e^{(\log N)^{1-p}}e^{C_{\ref{4.6}}t}
 X_0^N(P_t^{N,*}\phi)\cr
 &&\quad\quad+C_{\ref{4.6}}e^{C_{\ref{4.6}}t}||\phi||_{1/2}
 (\log N)^{(1-p)/2}X_0^N(1).
 \eeqlb
 \end{proposition}
 \begin{proposition}
 \label{prop4.5}
 For $p\geq 3$ there is a constant $C_{\ref{4.7}}(p)$ such that for
 all $\phi:\mbb R\rar\mbb R^+$, if $\ez=(\log N)^{-p}$, then
 \beqlb
 \label{4.7}
 E(X_{\ez}^N(\log N\phi f_0^N(\cdot,\xi_{\ez}^N)))\leq
 C_{\ref{4.7}}X_0^N(1)||\phi||_{1/2}(\log
 N)^{(1-p)/2}+C_{\ref{4.7}}X_0^N(\phi).
 \eeqlb
 \end{proposition}
Let $\sup_{K,T}$ indicate a supremum over all $X_0^N\in M(\mbb
S_N)$, $\phi:\mbb R\rar\mbb R$ and $t\geq0$ satisfying $X_0^N(1)\leq
K$, $||\phi||_{\textrm{Lip}}\leq K$ and $t\leq T$.
\begin{remark}
Note that if $||\phi||_{\textrm{Lip}}\leq K$, then
$||\phi||_{\underline{\alpha}}\leq 2K$ for any
$0<\underline{\alpha}<1$.
\end{remark}
\begin{proposition}
\label{prop4.7} For every $K,T>0$ and $0<p<2$,
 \bgeqn
\label{4.7a}\lim_{N\rar\infty}\sup_{K,T}E\left(\left|\int_0^tX_s^N(\log
N\phi^2f_0^N(\cdot,\xi_s^N))-p_1(0)^{-1}X_s^N(\phi^2)\right|^p\right)=0
\edeqn and for i=2,\,3, \bgeqn \label{4.7b}
\lim_{N\rar\infty}\sup_{K,T} E\left(\left|D_t^{N,i}
-\int_0^t\theta_{i-2}\gamma^*X_s^N(\phi)ds\right|^p\right)=0. \edeqn
\end{proposition}

Recall the rescaled Lotka-Volterra models in Section \ref{ourmodel}
and  assume (A2) holds. Also recall the 1-biased voter model and
0-biased voter model with rates $\nu=\nu_N$ and $b=b_N$ defined in
the last subsection. Set $\bar{\xi}_t^N(x)=\bar{\xi}_t(Nx)$ and
$\underline{\xi}_t^N(x)=\underline{\xi}_t(Nx)$ for $x\in \mbb S_N$.
Thus the rate function of $\bar{\xi}_t^N$ is given by\beqnn
 \bar{c}(x,\xi)=\begin{cases}
 (\nu_N+b_N)f_1^N(x,\xi)& \textrm{ if }\xi(x)=0,\cr
 \nu_N f_0^N(x,\xi)& \textrm{ if }\xi(x)=1
                      \end{cases}
 \eeqnn
and the rate function of $\underline{\xi}_t^N(x)$ is given by \beqnn
 \underline{c}(x,\xi)=\begin{cases}
 \nu_N f_1^N(x,\xi)& \textrm{ if }\xi(x)=0,\cr
 (\nu_N+b_N) f_0^N(x,\xi)& \textrm{ if }\xi(x)=1.
                      \end{cases}
 \eeqnn
Assume $N$ is large enough ($N\geq N_0$) so that $\nu_N>0$ and
$b_N>1$. As in the last subsection, we may construct the three
processes on one probability space so that
$\underline{\xi}^N_0=\hat{\xi}^N_0=\bar{\xi}^N_0$ and \bgeqn
\label{6.1}\underline{\xi}^N_t\leq\hat{\xi}^N_t\leq\bar{\xi}^N_t\textrm{
for all }t\geq0.\edeqn Define
$$\bar{X}_t^N=\frac{1}{N'}\sum_{x\in\mbb S_N}\bar{\xi}_t^N(x)\dz_{x}
\textrm{ and }\underline{X}_t^N=\frac{1}{N'}\sum_{x\in\mbb
S_N}\underline{\xi}_t^N(x)\dz_{x}.
$$
It follows that
 \bgeqn
 \label{6.2}
 \underline{X}_t^N\leq X_t^N\leq\bar{X}_t^N\textrm{ for all }t\geq0.
 \edeqn

Keep Remark \ref{remark5.3} in mind. Applying Proposition
\ref{prop5.4} gives that there are constants $C_{\ref{6.3}}$ and
$C_{\ref{4.1}}$ such that for all $N\geq N_0$ and $t\geq0$,
 \beqlb
 \label{6.3}
 E(\bar{X}_t^N(1))&\leq&
 C_{\ref{6.3}}e^{C_{\ref{6.3}}t}\bar{X}_0^N(1),\\
 \label{6.4}
 E(\bar{X}_t^N(1)^2)&\leq&
 C_{\ref{6.3}}e^{C_{\ref{6.3}}t}(\bar{X}_0^N(1)^2
 +\bar{X}_0^N(1))
 \eeqlb
 and if $g$ is as in (\ref{4.1}), then
 \beqlb
 \label{6.5}
 (\log N)E(\bar{X}_t^N(f_0^N(\cdot,\bar{\xi}_t^N)))
 &\leq& g(t)X_0^N(1),\\
 \label{6.6}
 (\log N)E(\bar{X}_t^N(1)\bar{X}_t^N(f_0^N(\cdot,\bar{\xi}_t^N)))&\leq&
 g(t)(X_0^N(1)^2+X_s^N(1)).
 \eeqlb
Typically, we have there exists a constant $C_{\ref{6.7}}$ such that
 \bgeqn
 \label{6.7}
 E(\bar{X}_t^N(1))-E(\underline{X}_t^N(1))\leq C_{\ref{6.7}}[(\log
 N)^{-2}+t]X_0^N(1),\quad0\leq t\leq1 \edeqn
whose counterpart in \cite{[CP08]}  is (6.7). We first prove
Proposition \ref{prop4.3}. In fact, we only give an outline.

{\textit{ Proof of Proposition \ref{prop4.3}}.} With inequalities
(\ref{6.3}), (\ref{6.4}) and the coupling (\ref{6.2}) in hand, part
(a) follows from the strong $L^2$ inequality for non-negative
submartingales and the fact that $\bar{X}_t^N(1)^2$ is a
submartingales; see Remark \ref{remarksec5} and (5.29) of
\cite{[CP08]}.   For part (b), if we have similar results to those
in Proposition 6.1 of \cite{[CP08]}, then part (b) follows from
Remark \ref{remark5.3}. But the proof of Proposition 6.1 of
\cite{[CP08]} works here if we replace (5.40) there by (\ref{5.40})
in the last subsection; see Remark \ref{remarksec5}. \qed

{\textit{ Proof of Proposition \ref{prop4.4}}}. Recall that
$\bar{\xi}_t$ is the biased voter model with rates
$\nu=N-\bar{\theta}\log N$ and $b=2\bar{\theta}\log N$, and
$\bar{\xi}_t^N(x)=\bar{\xi}_t(Nx)$, $x\in\mbb S_N$. For $\psi:\mbb
R\rar\mbb R^+$, define $\phi:\mbb Z\rar\mbb R^+$ by
$\phi(x)=\psi(x/N)$. Then $||\phi||_{\infty}=||\psi||_{\infty}$ and
for $0<\underline{\alpha}<1$, \beqnn\sup_{x\neq
y,|x-y|\leq1}\frac{|\phi(x)-\phi(y)|}{|x-y|^{\underline{\alpha}/2}}&\leq&
\sup_{x\neq y,|x-y|\leq1}\frac{|\phi(x)-\phi(y)|}{|x-y|^{1/2}}\cr
&\leq& N^{-1/2}\sup_{x\neq
y,|x-y|\leq1/N}\frac{|\psi(x)-\psi(y)|}{|x-y|^{1/2}}.\eeqnn Thus
$||\phi||_{\underline{\alpha}/2}\leq N^{-1/2}||\psi||_{1/2}.$ Note
that ${\cal A}^*_N\psi(x)=(N+\bar{\theta}\log N)\sum_{y\in\mbb
S_N}p_N(y-x)\psi(y)$ with semigroup $P_t^{N,*}$ and ${\cal
A}^*\phi(x)=(N+\bar{\theta}\log N)\sum_yp(y-x)\phi(y)$ with
semigroup $P_t^*$; see (\ref{defA*}) for the definition of ${\cal
A}^*$. We have that $P_t^*\phi(x)=P_t^{N,*}\psi(x/N)$ and
$\bar{\xi}^N_t(\psi)=\bar{\xi}_t(\phi)$. According to (\ref{5.23}),
we obtain
 \beqnn
 E(\bar{\xi}^N_t(\psi))\leq e^{b^{1-p}+(1+\kappa_p)t}
 \left(\bar{\xi}^N_0(P_t^{N,*}(\psi))+\left[\kappa_pb^{2-p}
 ||\psi||_{\infty}+B_pN^{-1/2}||\psi||_{1/2}
 \right]|\bar{\xi}^N_0| \right).
 \eeqnn
Since $p\geq3$, Remark \ref{remark5.3} implies
$\kappa_pb^{2-p}+B_pN^{-1/2}=O((\log N)^{(1-p)/2})$ as
$N\rar\infty.$ Then the fact that $\bar{\theta}\geq1$ implies
$b\geq\log N$  and the coupling (\ref{6.2}) yield the desired
inequality (\ref{4.6}). \qed

{\textit{ Proof of Proposition \ref{prop4.5}}}. Let $\ez=b^{-p}$.
According to Remark \ref{remarksec5}, we may use (5.32) of
\cite{[CP08]} to obtain that
 $$
 E(X_{\ez}^N(b\phi f_0(\xi_{\ez}^N)))\leq
  E(\bar{X}_{\ez}^N(b\phi
  f_0(\bar{\xi}_{\ez}^N)))+2b||\phi||_{\infty}
  (E(\bar{X}_{\ez}^N(1)-X_{\ez}^N(1))).
 $$
Applying (5.62) of \cite{[CP08]} and (\ref{5.40}) gives
$$
 E(X_{\ez}^N(b\phi f_0(\xi_{\ez}^N)))\leq
  (6eb^{2-p}||\phi||_{\infty}+B_p N^{-1/2}||\phi||_{1/2})X_0^N(1)
  +\kappa_p X_0^N(\phi).
 $$
Then Remark \ref{remark5.3} yields (\ref{4.7}). \qed

 We will give
the proof of Proposition \ref{prop4.7} in the final subsection. In
the next subsection with the help of the four propositions in this
subsection we prove Theorem \ref{main2}.

\subsection{Convergence Theorem}
In this subsection, we follow the strategy in the Section 4 of
\cite{[CP08]} to obtain Theorem \ref{main2}. First, we check the
compact containment condition.
\begin{proposition}
\label{Prop4.12} For all $\ez>0$ there is an $n\in \mbb N$, so that
$$
\sup_NP\left(\sup_{t\leq\ez^{-1}}X_t^N(B(0,n)^c)>\ez\right)<\ez.
$$
\end{proposition}
{\bf Proof. }The proof is similar to that for Proposition 4.12 of
\cite{[CP08]}. We only give an outline here. Recall that $b(N)=N$.
Let $h_n:\mbb R^d\rar[0,1]$ be a $C^{\infty}$ function such that
 $$
 1_{\{|x|>n+1\}}\leq h_n(x)\leq 1_{\{|x|>n\}}
 $$
and
 $$
 \sup_n\sum_{i,j,k\leq d}||(h_n)_i||_{\infty}+||(h_n)_{ij}||_{\infty}
 +||(h_n)_{ijk}||_{\infty}\equiv C_h<\infty. $$
By the semimartingale decomposition
$$
\sup_{t\leq T}X_t^N(h_n)\leq X_0^N(h_n)+\sum_{i=1}^3 \sup_{t\leq
T}|D_t^{N,i}(h_n)|+\sup_{t\leq T}|M_t^N(h_n)|.
$$
We need to check the right hand side tends to zero as
$N,n\rar\infty$. Let
$$
\eta_N:=\sup_n||{\cal
A}_N(h_n)-\frac{\sigma^2\Delta^{1/2}h_n}{2}||_{\infty}.
$$
Then $\lim_{N\rar\infty}\eta_N=0$ by Lemma \ref{conge}. Note that
 \beqnn
&&\frac{1}{N'}\sum_{x,y}|h_n(x)-h_n(y)|p_N(x-y)\xi_s^N(y)\\
&&\qquad\leq
\frac{||h_n||_{\underline{\alpha}}}{N'}\sum_{y}\sum_{x}|x-y|^{\underline{\alpha}}p_N(x-y)\xi_s^N(y)\\
&&\qquad\leq\frac{C_h|p|_{\underline{\alpha}}}{N^{\underline{\alpha}}}X_s^N(1).
 \eeqnn
Set $\eta'_N(T)=C_{\ref{4.2}}(T)(\eta_N+\bar{\theta}C_h\log
N|p|_{\underline{\alpha}}/{N^{\underline{\alpha}}})T )$. We have, as
in the deviation of (4.17) in \cite{[CP08]}
 \beqlb
 \label{4.16}
&& E\left(\sup_{t\leq T}X_t^N(h_n)\right)\leq X_0^N(h_n)+2(\la
M^N(h_n)\ra_T)^{1/2}+\eta'_N X_0^N(1)\cr
&&\quad+C_h\int_0^TE(X_s^N(h_{n-1}))ds+2\bar{\theta}\int_0^TE(X_s^N(h_n\log
Nf_0^N(\xi_s^N)))ds.
 \eeqlb
Applying Proposition \ref{prop4.5} and (\ref{4.2}), we obtain the
last integral above is bounded by
 \bgeqn
 \label{4.17}
 \eta''_N(T)X_0^N(1)+C_{\ref{4.7}}\int_0^TE(X_s^N(h_n))ds,
 \edeqn
where $\eta''_N(T)=C_{\ref{4.2}}(T)[(\log N)^{-2}+C_{\ref{4.7}}C_h
T/\log N].$ By Lemma \ref{estimateMNt} and (\ref{4.2}) there is a
constant $C_{\ref{4.11}}(T)$ such that if $\phi_s=\psi$, then for
any $\underline{\alpha}<1$ and $0\leq s\leq T,$
 \bgeqn
 \label{4.11}
 E(|m_{1,s}^N|+|m_{2,s}^N|)\leq C_{\ref{4.11}}(T)||\phi||^2_{\underline{\alpha}}
 (\log N/N^{\underline{\alpha}})X_0^N(1).
 \edeqn
Then the above inequality, (\ref{4.17}) and  Lemma \ref{estimateMNt}
gives (recall $N/N'=\log N$)
 \bgeqn \label{4.18}
 E(\la M^N(h_n)\ra_T)\leq
 \eta_N'''(T)X_0^N(1)+2C_{\ref{4.7}}\int_0^TE(X_s^N(h_n))ds,
 \edeqn
where $\eta_N'''(T)=2\eta_N''(T)+C_{\ref{4.11}}(T)TC_h^2\log
N/N^{\underline{\alpha}}.$ Finally, let $B_t^{N,*}$ be the
continuous random walk with semigroup $P_t^{N,*}$ defined before
Proposition \ref{prop4.4}, $B_0^{N,*}=0$. Note that
$$
P\left(|B_s^{N,*}|\geq
\frac{n-1}{2}\right)=P\left(|B^0_{(N+\bar{\theta}\log N)s}|\geq
\frac{N(n-1)}{2}\right).
$$
Since $b(t)=l(t)=t$, Proposition \ref{propincrease} yields that the
left hand side above goes to 0 uniformly in $N\in\mbb N$ and $0\leq
s\leq T$ as $n\rar\infty$. Thus with the help of Proposition
\ref{prop4.4} and the inequalities (\ref{4.16}), (\ref{4.17}),
(\ref{4.18}) we can conclude: for any $T,\ez>0$ there is an $N_0$
such that
$$
\textrm{for } N\geq N_0, n\geq N_0, E(\sup_{t\leq T}X_t^N(h_n))<\ez.
$$
The desired result is immediate. \qed

{\textit{ Proof of Theorem \ref{main2}. }} In fact, we have already
completed all tasks. First, with (\ref{4.4}) and (\ref{4.5}) in
hand, by the same arguments as those in the proof of Lemma 4.10 of
\cite{[CP08]}, we have there exists a constant $C_{\ref{4.12}}(T)$
such that for all $0\leq s\leq t\leq T$,
 \bgeqn
 \label{4.12}
E\left(\left[\int_s^t X_r^N(\log N f_0^N(\xi_r^N))dr\right]^2\right)
\leq C_{\ref{4.12}}(T)(t-s)^{4/3}(X_0^N(1)^2+X_0^N(1)).
 \edeqn
Now, recall the decomposition of $X_t^N(\phi_t)$ in Section
\ref{secsemi}. With the help of Lemma \ref{estimateMNt} and
(\ref{4.12}), by the the same  arguments as those in the proof of
Proposition 4.11 of \cite{[CP08]}, for each $\phi\in C_b^{1,3}({\mbb
R_+}\times\mbb R)$, each of families $\{X_{\cdot}^N(\phi),N\in\mbb
N\}$, $\{D_{\cdot}^{N,i},N\in\mbb N\}$, $i=1,2,3,$ $\{\la
M^N(\phi\ra)_{\cdot},N\in\mbb N\}$, and
$\{M_{\cdot}^N(\phi),N\in\mbb N\}$ is C-tight in $D([0,\infty),\mbb
R)$. The C-tightness of $\{P_N,N\in\mbb N\}$ is now immediate from
Proposition \ref{Prop4.12} and Theorem II.4.1 of \cite{[P02]}. Then
to check any limit point of $\{P_N\}$ is the law claimed in the
Theorem,  one
 can follow the same arguments as those in the proof of proposition
 4.2 of \cite{[CP08]}, using Proposition \ref{prop4.7} above.  \qed

\subsection{Proof of Proposition \ref{prop4.7}}
For $N$ fixed, let $\hat{\xi}_t$ be the rate
$\nu_N=N-\bar{\theta}\log N$ voter model on $\mbb Z$ with rate as in
(\ref{5.1}) for $b=0$ and $\nu=\nu_N$. Define
$\hat{\xi}_t^N(x)=\hat{\xi}_t(xN)$, $x\in\mbb S_N$, the rate $\nu_N$
voter model on $\mbb S_N$.  Recall the independent and coalescing
random walks system $\{B_t^x\}$ and $\{\hat{B}_t^x\}$ defined in
Section \ref{MR}. We need to introduce their rescaled versions as
follows: for $x,y\in\mbb S_N$,
 \bgeqn
 \label{7.13}
 B_t^{N,x}=B^{xN}_{\nu_Nt}/N,\quad
 \hat{B}_t^{N,x}=B_{\nu_Nt}^{xN}/N,
 \edeqn
and
 $$
 \tau^N(x,y)=\tau(Nx,Ny)/\nu_N,\quad\hat{\tau}^N(x,y)=\hat{\tau}(Nx,Ny)/\nu_N.
 $$
Define
 $$
 \varepsilon(t)=\sup_{x\in\mbb Z}|tp_t(0,x)-p_1(x/t)|\vee(1/t^2).
 $$
By Proposition \ref{CPtostable} $\varepsilon(t)\rar0$ as
$t\rar\infty$.
 Then for each $k\in\mbb Z^+$, there exists a
$t(k)$ such that for $t>t(k)$, $\varepsilon(t)\leq 1/k$. Define
\beqlb\varepsilon'(t)=\begin{cases}1,&0\leq t\leq t(1),\\
1/k,& t(k)<t\leq t(k+1).
\end{cases}\eeqlb
Then $\varepsilon'(t)\downarrow0$ as $t\rar\infty$ and
$\varepsilon'(t)\geq\varepsilon(t)$ for $t>t(1)$.  Let
$\hat{\eta}_N=e^{-\sqrt{\log N}}$ and
$a_N=\nu_N(2-\hat{\eta}_N)/\log N$ and $$ \ez_N'=({\log\log
N})^{-1}\vee\sqrt{\varepsilon'(a_N/\log\log N)}.
$$
Then \beqnn&&\ez_N:=\left(\varepsilon(a_N\ez'_N)/\ez'_N
 +\frac{\log\log N}{\log N}\right)\\&&\quad\leq\varepsilon'(a_N\ez'_N)
\left(\sqrt{\varepsilon'(a_N/\log\log N)}\right)^{-1}+\frac{\log\log N}{\log N}\\
&&\quad\leq\sqrt{\varepsilon'(a_N/\log\log N)}+\frac{\log\log
N}{\log N} \rar0\eeqnn as $N\rar\infty$. Define the sequences
 \bgeqn
 \label{7.3}
 t_N=\frac{\ez'_N}{\log N},\quad K_N=(\log
 N)^{1/2},\quad\dz_N=K_Nt_N.
 \edeqn
We assume that $N$ is large enough so that $\ez'_N\vee t_N\vee
\dz_N\leq 1$ and $\dz_N/\ez'_N\rar0$ as $N\rar\infty$. The following
lemma is a version of Lemma 7.6 of \cite{[CP08]}.
 \begin{lemma}
 \label{lemma7.6}
There is a constant $C_{\ref{7.14}}$ such that
 \beqlb
 \label{7.14}
 &&\frac{\log N}{N'}\sum_{x,e}p_N(e)P\left(\hat{\xi}_0^N(B_{t_N}^{N,x})=
 \hat{\xi}_0^N(B_{t_N}^{N,x+e})=1,\tau^N(x,x+e)>t_N\right)\cr
 &&\quad\quad\leq C_{\ref{7.14}}(\ez'_N)^{-1}\int\int_{|w-z|\leq\dz_N}d
 \hat{X}^N_0(w)d\hat{X}^N_0(z)+C_{\ref{7.14}}\ez_N
 \hat{X}_0^N(1)^2.
 \eeqlb\end{lemma}
{\bf Proof. }By translation invariance and symmetry, the left side
of (\ref{7.14}) is
 \beqlb
 \label{7.15}
 &&(N')^{-2}\sum_{w,z}\hat{\xi}_0^N(w)\hat{\xi}_0^N(z)\sum_ep_N(e)\cr
 &&\qquad\quad\times\left[\sum_x NP(B_{t_N}^{N,0}=w-x,
 B_{t_N}^{N,e}=z-x,
 \tau^N(0,e)>t_N)\right]\cr
 &&\quad=(N')^{-2}\sum_{w,z}\hat{\xi}_0^N(w)\hat{\xi}_0^N(z)
 \sum_ep_N(e)NP(B_{2t_N}^{N,e}=z-w,\tau_0^{N,e}>2t_N)\cr
 &&\quad\equiv\Sigma_d^N+\Sigma_c^N,
 \eeqlb
where $\tau_0^{N,e}=\inf\{s:B_s^{N,e}=0\}$, and $\Sigma_d^N$,
respectively, $\Sigma_c^N,$ denotes the contribution to (\ref{7.15})
from $w,z$ satisfying $|w-z|\leq K_Nt_N$, respectively,
$|w-z|>K_Nt_N$. Let
$$
\tilde{P}((B_{\cdot}^N,\tau_0^N)\in\cdot)=\sum_ep_N(e)P((B^{N,e}_{\cdot},
\tau_0^{N,e})\in\cdot).
$$
For $\Sigma_d^N$, use (\ref{boundtransition}) and the Markov
property at time $t_N$ to see that
 \beqnn
 && N\tilde{P}^N(B_{2t_N}^N=z-w,\tau_0^N>2t_N)\\
 &&\quad\quad\leq N\tilde{E}(P(B_{t_N}^{N,0}=z-w-B_{t_N}^N(w));
 \tau_0^N>t_N)\\
 &&\quad\quad\leq CN\tilde{P}(\tau_0^N>t_N)(\nu_Nt_N)^{-1}\\
 &&\quad\quad\leq C\frac{NH(\nu_Nt_N)}{\nu_Nt_N}.
  \eeqnn
By (\ref{2.4}), there is a constant $C_{\ref{7.16}}$ such that
 \bgeqn
 \label{7.16}
 \Sigma_d^N\leq C_{\ref{7.16}}(\ez'_N)^{-1}\int\int_{|w-z|\leq K_Nt_N}
 d\hat{X}^N_0(w)d\hat{X}^N_0(z).
 \edeqn
It is more complicated to bound $\Sigma_c^N$.  Using the Markov
property at time $\hat{\eta}_Nt_N$ gives
 \beqnn
 &&\tilde{P}^N\left(B_{2 t_N}^N=w-z,\tau_0^N>2t_N\right)\\
 &&\qquad\leq\tilde{P}\left(\tau_0^N>\hat{\eta}_Nt_N,
 |B_{\hat{\eta}_Nt_N}^N|>\frac{K_Nt_N}{2}\right)\sup_{x'}
 P\left(B_{(2-\hat{\eta}_N)t_N}^{N,0}=x'\right)\cr
 &&\qquad\quad+\tilde{P}\left(
 P\left(B_{(2-\hat{\eta}_N)t_N}^{N,0}=w-z-B_{\hat{\eta}_Nt_N}^N\right);
 \tau_0^N>\hat{\eta}_Nt_N,
 |B_{\hat{\eta}_Nt_N}^N|\leq\frac{K_Nt_N}{2}\right)\cr
 &&\qquad=\Sigma_{1c}^N+\Sigma_{2c}^N,\quad\qquad\textrm{say}.
 \eeqnn
Note that
 $$\tilde{P}\left(|B_{\hat{\eta}_Nt_N}^N|>\frac{K_Nt_N}{2}\right)=
 \sum_ep_N(e)P\left(|B_{N\hat{\eta}_Nt_N}^0+e|>\frac{NK_Nt_N}{2}\right)$$
which is bounded by
$$
\frac{2|p|_{1/2}}{(NK_Nt_N)^{1/2}}
+P\left(|B_{N\hat{\eta}_Nt_N}^0|>\frac{NK_Nt_N}{4}\right).
$$
By Proposition \ref{propincrease},
$$
P\left(|B_{N\hat{\eta}_Nt_N}^0|>\frac{NK_Nt_N}{4}\right)\leq
\frac{4C_{\ref{increasein}}N\hat{\eta}_Nt_N}{{NK_Nt_N}}=
4C_{\ref{increasein}}\hat{\eta}_N/K_N.
$$
(Note that $l(t)=b(t)=t$.) Thus by (\ref{7.6})
 \beqlb\label{1c}
\Sigma_{1c}^N\leq
\frac{C\left(\hat{\eta}_N/K_N+1/(NK_Nt_N)^{1/2}\right)}{\nu_N(2-\hat{\eta_N})t_N}.
 \eeqlb
Let us consider $\Sigma_{2c}^N$. By the definition of
$\varepsilon(t)$ and  (\ref{estistable}) (recall
$d=\alpha=1$),\beqlb\label{inequ4.7}
 p_t(0,x)&\leq& \frac{\varepsilon(t)}{t}+\frac{p_1(x/t)}{t}\cr
 &\leq& \frac{1}{t}\left({\varepsilon(t)}+c_2\left(1\wedge
 \left|\frac{t}{x}\right|^2\right)\right).
\eeqlb Note that for $|w-z|>K_Nt_N$, on
$\left\{|B_{\hat{\eta}_Nt_N}^N|\leq\frac{K_Nt_N}{2}\right\}$,
$$|w-z-B_{\hat{\eta}_Nt_N}^N|^{-1}\leq\frac{2}{K_Nt_N}. $$
Thus by inequality (\ref{inequ4.7}), $\Sigma_{2c}^N$ is less than
$$
\left(\varepsilon(\nu_N(2-\hat{\eta}_N)t_N)
+c_2\left(1\wedge\left(\frac{2\nu_N(2-\hat{\eta}_N)}{NK_N}\right)^2\right)\right)
\frac{H(\nu_N\hat{\eta}_Nt_N)}{\nu_N(2-\hat{\eta}_N)t_N}.
$$
Thus by $a_N\ez'_N=\nu_N(2-\hat{\eta}_N)t_N$ and (\ref{2.4}),
 \beqlb\label{2c}
\Sigma_{2c}^N&\leq& C\left(\varepsilon(a_N\ez'_N)+1/K_N^2\right)
\frac{\log N}{\nu_N\ez'_N\log (\nu_N\hat{\eta}_Nt_N)}\cr &\leq&
 C\left(\varepsilon(a_N\ez'_N)/(N\ez_N')+(N\log N\ez_N')^{-1}\right)\cr
 &\leq &C\left(\varepsilon(a_N\ez'_N)/(N\ez_N')+\frac{\log\log N}{N\log
 N}\right)\cr
 &=& C\ez_N/N,
 \eeqlb
where $C$  may change its values from line to line and the second
inequality follows from
$$
\log(\nu_n\hat{\eta}_Nt_N)=\log (\ez_N')+\log (\nu_N)-\log\log
N-\sqrt{\log N}
$$  and $\lim_{N\rar\infty}\frac{N}{\nu_N}=1.$ With (\ref{7.16}),
(\ref{1c}) and (\ref{2c}) in hand,  (\ref{7.15}) yields the desired
result, (\ref{7.14}). \qed

For $\phi:\mbb R^2\rar\mbb R$, $\zeta\in\{0,1\}^{\mbb S_N}$ and
$X(\phi)=(1/N')\sum_x\phi(x)\zeta(x)$, define
 \beqnn
 \Delta_1^{N,+}(\phi,\zeta)&=& X(\log N\phi^2f_0^N(\cdot,\zeta))\\
 \Delta_2^{N,+}(\phi,\zeta)&=&
 \frac{1}{N'}\sum_x(1-\zeta(x))\phi(x)\log N f_1^N(x,\zeta)^2\\
 \Delta_3^{N,+}(\phi,\zeta)&=& X(\log N\phi
 f_0^N(\cdot,\zeta)^2)
 \eeqnn
and
 $$
 \Delta_j^N(\phi,\zeta)=\Delta_j^{N,+}(\phi,\zeta)\gamma_j
 X(\phi),\quad j=1,2,3,
 $$
where $\gamma_1=p_1(0)^{-1}$ and $\gamma_2=\gamma_3=\gamma^*$.
Define
 $$
 m(1)=2\quad\textrm{and}\quad m(2)=m(3)=1.
 $$
The following proposition is a version of Proposition 7.5 of
\cite{[CP08]}.
 \begin{proposition}
 \label{prop7.5}
 There is a constant $C_{\ref{7.12}}$ and a sequence
 $\eta_{\ref{7.12}}(N)\downarrow0$ such that for $j=1,2,3,$ if
 $\phi:\mbb R^2\rar\mbb R,$ then for any $0<\underline{\alpha}<1$
  \beqlb
  \label{7.12}
  |E(\Delta_j^N(\phi,\hat{\xi}_{t_N}^N))|
 &\leq&\eta_{\ref{7.12}}(N)\left(\hat{X}_0^N(1)+\hat{X}_0^N(1)^2\right)
 ||\phi||_{\underline{\alpha}}^{m(j)}\cr
 &&\quad+\frac{C_{\ref{7.12}}||\phi||_{\infty}^{m(j)}}{\ez_N'}
 \int\int_{|w-z|\leq\ez_N}d\hat{X}_0^N(w)d\hat{X}_0^N(z).
  \eeqlb
  \end{proposition}
{\bf Proof. } To prove the proposition, we can define
$\Sigma_j^{i,N}$, $i=1,2$ for $j=1$ and $i=1,2,3$ for $j=2,3$ as in
(7.20), (7.21) and (7.22) of \cite{[CP08]} and decompose each
$E(\Delta_j^{N,+})$ into a sum of those terms. We omit the
definitions and decompositions here, since they are the same. By
Lemma \ref{lemma7.6}, we can show that
 \beqlb
 \label{7.23}
 &&\Sigma_j^{2,N}\leq C_{\ref{7.14}}||\phi||_{\infty}^{m(j)}
  \left[(\ez'_N)^{-1}\int\int_{|w-z|\leq\dz_N}d
 \hat{X}^N_0(w)d\hat{X}^N_0(z)+\ez_N\hat{X}_0^N(1)^2\right].
 \eeqlb
For $\Sigma_j^{3,N}$, $j=2,3$, with Proposition \ref{prop2.2} in
hand, one can check that a similar conclusion to that in Lemma 2.5
of \cite{[CP08]} is available. Following the proof of Proposition
7.5 of \cite{[CP08]}, we have there exists a constant
$C_{\ref{7.24}}$ depending on $p(\cdot)$,
 \bgeqn
 \label{7.24}\Sigma_2^{3,N}+\Sigma_3^{3,N}\leq C_{\ref{7.24}}
||\phi||_{\infty}\hat{X}_0^N(1)(\log N)^{-1/2}.
 \edeqn
Now, we need to establish that there is a sequence $\eta(N)\rar0$
such that for $j=1,2,3$,
 \bgeqn
 \label{7.19}
 |\Sigma_j^{1,N}-\gamma_j\hat{X}_0^N(\phi)|\leq\eta(N)||\phi||_{\underline{\alpha}}^{m(j)}
 \hat{X}_0^N(1).
 \edeqn
 Let $e$ denote independent
random variable with  law $p(\cdot)$. First, $$P\left(
{B}_{t_N}^{N,e}>\sqrt{\ez_N'}\right)=P\left(
|{B}_{\nu_Nt_N}^{0}+e|>N\sqrt{\ez_N'}\right).$$ We also have
 \beqlb
 \label{bound1}
P\left( |{B}_{\nu_Nt_N}^{0}+e|>N\sqrt{\ez_N'}\right)&\leq&
\frac{2|p|_{\underline{\alpha}}}{(N\sqrt{\ez_N'})^{\underline{\alpha}}}+P\left(
|{B}_{\nu_Nt_N}^{0}|>N\sqrt{\ez_N'}/2\right)\cr
&\leq&\frac{2|p|_{\underline{\alpha}}}{(N\sqrt{\ez_N'})^{\underline{\alpha}}}+
\frac{C_{\ref{increasein}}\nu_Nt_N}{N\sqrt{\ez_N'}},
 \eeqlb
where  the second inequality follows from  Proposition
\ref{propincrease}. Typically, we have
 \bgeqn
 \label{bound2}
 P\left({B}_{t_N}^{N,0}>\sqrt{\ez_N'}\right)\leq \frac
 {C_{\ref{increasein}}\nu_Nt_N}{N\sqrt{\ez_N'}}
 = \frac{C_{\ref{increasein}}\nu_N\sqrt{\ez_N'}}{N\log N}.
 \edeqn
 Now, we
consider the case of $j=2$. By the same arguments as in
\cite{[CP08]}, we can show
 \beqnn
 &&|\Sigma_2^{1,N}-\gamma^*\hat{X}_0^N(\phi)|\\
 &&\leq \frac{1}{N'}\sum_w\hat{\xi}_0^N(w)\log N
 E\bigg{(}|\phi(w-\hat{B}_{t_N}^{N,e})-\phi(w)|;
 \hat{\tau}^N(0,e)\wedge\hat{\tau}^N(0,f)>t_N,\\
 &&\qquad\quad\hat{\tau}^N(e,f)\leq
 t_N\bigg{)}+\left|\frac{1}{N'}\sum_w\hat{\xi}_0^N(w)\phi(w)
 (q_{\nu_Nt_N}\log N -\gamma^*)\right|\\
 &&\leq ||\phi||_{\underline{\alpha}}\hat{X}_0^N(1)\log
 N\left(\sqrt{\ez_N'}\right)^{\underline{\alpha}}q_{\nu_Nt_N}
 +2||\phi||_{\infty}\hat{X}_0^N(1)\log N
 P\left(|{B}_{t_N}^{N,e}|>\sqrt{\ez_N'}\right)^{1/2}q_{\nu_Nt_N}^{1/2}\\
 &&\qquad+||\phi||_{\infty}\hat{X}_0^N(1)|\log N
 q_{\nu_Nt_N}-\gamma^*|,
 \eeqnn
 where the second inequality follows from Cauchy-Schwarz
 inequality and considering the cases $|{B}_{t_N}^{N,e}|>\sqrt{\ez_N'}$
  and $|{B}_{t_N}^{N,e}|\leq \sqrt{\ez_N'}$. Thus by (\ref{2.4}),
  (\ref{bound1}) and Proposition \ref{Prop2.1}, there exists a
  sequence $\eta_{\ref{7.26}}(N)$ which goes to 0 as $N\rar\infty$
  such that
  \bgeqn
  \label{7.26}
  |\Sigma_2^{1,N}-\gamma^*\hat{X}_0^N(\phi)|\leq
  \eta_{\ref{7.26}}(N)||\phi||_{\underline{\alpha}}\hat{X}_0^N(1).
  \edeqn
By replacing $\hat{B}_{t_N}^{N,e},{B}_{t_N}^{N,e}$ with
$\hat{B}_{t_N}^{N,0},{B}_{t_N}^{N,0}$ respectively,  the same
argument as that above gives the same bound for
$|\Sigma_3^{1,N}-\gamma^*\hat{X}_0^N(\phi)|$. Typically, inequality
$(\ref{bound1})$ could be simplified. Next, we turn to
$\Sigma_2^{1,N}$. Following the strategy of the proof for term on
$\Sigma_2^{1,N}$, we have that
 \beqnn
 &&|\Sigma_1^{1,N}-p_1(0)^{-1}\hat{X}_0^N(\phi^2)|\\
 &&\quad=\left|\frac{1}{N'}\sum_w\hat{\xi}_0^N(w)\left[\log
 NE\left(\phi^2(w-B_{t_N}^{N,0});\tau^N(0,e)>t_N\right)
 -p_1(0)^{-1}\phi^2(w)\right]\right|\\
 &&\quad\leq \frac{1}{N'}\sum_w\hat{\xi}_0^N(w)\left[\log
 NE\left(\left|\phi^2(w-B_{t_N}^{N,0})-\phi^2(w)\right|;\tau^N(0,e)>t_N\right)
\right]\\
&&\qquad+\frac{1}{N'}\sum_w\hat{\xi}_0^N(w)\phi^2(w)|\log
NP(\tau^N(0,e)>t_N)-p_1(0)^{-1}|\\
&&\quad\leq \bigg{(} 2||\phi||_{\underline{\alpha}}\log
 N\left(\sqrt{\ez_N'}\right)^{\underline{\alpha}}H({\nu_Nt_N})
 +2||\phi||_{\infty}\log N
 P\left(|{B}_{t_N}^{N,0}|>\sqrt{\ez_N'}\right)^{1/2}H({\nu_Nt_N})^{1/2}\\
 &&\qquad+||\phi||_{\infty}
 |\log NH({\nu_Nt_N})-p_1(0)^{-1}|\bigg{)}||\phi||_{\infty}\hat{X}_0^N(1).
 \eeqnn
According to (\ref{bound1}) and (\ref{2.4}), we can conclude
 \bgeqn
 \label{7.27}
 |\Sigma_1^{1,N}-p_1(0)^{-1}\hat{X}_0^N(\phi^2)|\leq
 \eta_{\ref{7.27}}\hat{X}_0^N(1)||\phi||_{\underline{\alpha}}^2,
 \edeqn
where $\eta_{\ref{7.27}}\rar0$ as $N\rar\infty$. Thus we get the
(\ref{7.19}). By decompositions in (7.18) of \cite{[CP08]}, we
obtain the desired result.\qed

 With Proposition \ref{prop7.5} in hand, Proposition \ref{prop4.7}
 follows from the following two propositions which are analogous to
 Proposition 7.1 and Proposition 7.2 in \cite{[CP08]} and a similar argument
 to that in Section 8 of \cite{[CP08]}.
 \begin{proposition}
 \label{prop7.1}
 There is a  constant $C_{\ref{7.4}}(K)$ and sequence
 $\eta_{\ref{7.4}}(N)\downarrow0$ such that for all $\phi:\mbb
 R\rar[ 0,\infty)$ satisfying $||\phi||_{\textrm{Lip}}\vee X_0^N(1)\leq
 K$ and $j=1,2,3$,
 \beqlb
 \label{7.4}
 |E(\Delta_j^N(\phi,\xi_{t_N}^N))|&\leq&
 C_{\ref{7.4}}(K)\bigg{(}\eta_{\ref{7.4}}(N)\left(X_0^N(1)+X_0^N(1)^2\right)\cr
 &&\qquad+(\ez_N')^{-1}\int\int_{|w-z|\leq\dz_N}dX_0^N(w)dX_0^N(z)\bigg{)}.
 \eeqlb
 \end{proposition}
{\bf Proof. }First, we can obtain follow the strategy in the proof
of Lemma 7.8 in \cite{[CP08]} to obtain an analogous result to that
in Lemma 7.8 of \cite{[CP08]}. Then with our coupling, (\ref{6.7})
and  Proposition \ref{prop7.5} in hand, following the argument in
\cite{[CP08]}, one can get the desired result.  \qed

\begin{proposition}
\label{prop7.2} There is a constant $C_{\ref{7.5}}$ such that for
all $0\leq t\leq T$,
 \beqlb
 \label{7.5}
 &&E\bigg{(}\int\int_{|w-z|\leq\dz_N}dX_0^N(w)dX_0^N(z)\bigg{)}\cr
 &&\qquad\leq
 C_{\ref{7.5}}e^{C_{\ref{7.5}}T}(X_0^N(1)+X_0^N(1)^2)\cr
 &&\qquad\qquad\times\left[\frac{\dz_N}{\dz_N+t}(1+t^{2/3})+\dz_Nt^{-1/3}
 \log (1+\frac{t}{\dz_N})
 \right].
 \eeqlb
\end{proposition}
The proof of Proposition \ref{prop7.2} is also exactly the same with
that of Proposition 7.2 of \cite{[CP08]}. In fact, we only need to
prove the following random walk estimate which is a version of
Corollary 7.9 of \cite{[CP08]} and can be deduced directly from
(\ref{7.6}) and Proposition \ref{lemma7.3}. Let $B_t^{N,*}$ be the
random walk with semigroup $(P_t^{N,*},t\geq0)$ from Proposition
\ref{prop4.4}, at rate $\nu_N+b_N=N+\bar{\theta}\log N$,
$B_{\cdot}^{N,*}$ takes steps with $p_N(\cdot)$ and $B_0^{N,*}=0$.
 \begin{corollary}
 \label{cor7.9}
 (a) For all $x\in \mbb S_N$ and $t\geq0$,
  \bgeqn
  \label{7.30}
  P(B_t^{N,*}=x)\leq \frac{C_{\ref{7.6}}}{1+Nt}.
  \edeqn
(b) Assume $\dz_N'\downarrow0$ and $N\dz_N'\rar\infty$. For each
$K>0$ there is a constant $C_{\ref{7.31}}(K)>0$ such that
 \bgeqn
 \label{7.31}
 \inf_{N\geq1,w\in\mbb S_N,|w|\leq
 K\dz_N'}N\dz_N'P(B_{2\dz_N'}^{N,*}=w)\geq C_{\ref{7.31}}(K)>0.
 \edeqn
\end{corollary}
Now, one follows the argument in \cite{[CP08]} to get Proposition
\ref{prop7.2}. To obtain Proposition \ref{prop4.7}, the following
arguments are similar to those in Section 8 of \cite{[CP08]}.  We
omit it here.

\section{Voter Model's Asymptotics}
In  this section, we will prove Theorem \ref{voterasy} and we assume
that
 assumption (A1) holds with $b(t)=t^{1/\alpha}$.
Recall that $p_t=P(|\xi_t^0|>0)$. Our first object is to prove that
\begin{alignat}{2}
\label{obj1} p_t&=O\left(\frac{\log t}{t}\right)\quad&\textrm{as
}~t\rar\infty\quad& d=\alpha,\cr
   &=O(t^{-1})\quad&\textrm{as }~t\rar\infty\quad& d>\alpha.
\end{alignat}
The asymptotics above are similar to the results in Theorem 1 of
\cite{[BG80]}. Note that Theorem 1 of \cite{[BG80]} could be proved
under the assumption that the underlying motion has finite variance
and one only need to modify the proof of Lemma 5 of \cite{[BG80]};
see Lemma 2 of \cite{[BCL01]}. For our purpose we also need to
generalize the  asymptotic results in (14) of \cite{[BG80]}.

\medskip

Recall that  $\{B_t^x,x\in\mbb Z^d\}$ is a collection of rate-one
independent stable random walks with $B_0^x=x$. Let
$p_t(x,y)=P(B_t^x=y)$ denote the transition function of $\{B_t^x\}$.
Define the mean range of the stable random walk $B_t^0$ by
$$
R(t)=E\left(\sum_x1_{\{B_s^0=x \textrm{ for some } s\leq
t\}}\right).$$ By the  results for the range of the discrete time
stable random walk in \cite{[LR91]}, we see
 \begin{alignat}{2}
 \label{Range}
\lim_{t\rar\infty}\frac{R(t)}{t/\log
t}&=p_1(0)^{-1}\qquad&d=\alpha,\cr
\lim_{t\rar\infty}\frac{R(t)}{t}&=\gamma_e\qquad&d>\alpha.
 \end{alignat}
With this in hand, one can generalize the asymptotics results in
(14) of \cite{[BG80]}.   Now, to prove (\ref{obj1}) we only need to
prove some analogous results to those in Lemma 5 of \cite{[BG80]}.
Set $G_t(x)=\int_0^tp_s(0,x)ds$ and let
$\tau(x)=\inf\{t\geq0:B_t^x=0\}$, define $H_t(x)={P}(\tau(x)\leq
t).$
\begin{lemma}
\label{aymgreen} If $x\in \mbb Z^d$ with $|x|=r$, then there is a
constant $C_{d,\alpha}>0$ such that
\begin{alignat*}{2}
H_{r^{\alpha}}(x)&\geq C_{d,\alpha}/\log r\qquad &d=\alpha,\\
          &\geq C_{d,\alpha} r^{\alpha-d}\qquad &~d>\alpha.
\end{alignat*}
\end{lemma}
{\bf Proof. } We first consider the asymptotics for the Green's
function. According to  (\ref{tranapp}) and (\ref{estistable}), when
$r$ large enough,
$$
G_{r^{\alpha}}(x)=\int_0^{r^{\alpha}}p_s(0,x)ds\geq
c_1\int_{r^{\alpha}/2}^{r^{\alpha}}\frac{s}{r^{d+\alpha}}
ds-\int_{r^{\alpha}/2}^{r^{\alpha}}s^{-d/\alpha}ds.
$$
A bit of calculation show that there exist a constant
$\bar{C}_{d,\alpha}>0$ such that
\begin{alignat*}{2}
G_{r^{\alpha}}(x)&\geq \bar{C}_{d,\alpha}r^{\alpha-d}\quad&d>\alpha,\\
                 &\geq \bar{C}_{d,\alpha}\quad&d=\alpha.
\end{alignat*}
By (\ref{boundtransition}), we see that there exist constants
$\underline{C}_{d,\alpha}>0$ such that
\begin{alignat*}{2}
G_{r^{\alpha}}(0)&\leq\underline{C}_{d,\alpha}\quad&d>\alpha,\\
&\leq\underline{C}_{d,\alpha}\log r\quad&d=\alpha.
\end{alignat*}
Then the desired result follows from inequality $ H_t(x)\geq
G_t(x)/G_t(0). $\qed

 Now, one can follow the arguments in Section 3 of \cite{[BG80]} to
 obtain (\ref{obj1}) (Note that when prove an analogous result to that
  in Lemma 4 of \cite{[BG80]} one may need to set
  $s_t=d[(2p_t^{-1})^{1/d}]^{\alpha}$.)
With (\ref{obj1}), Theorem \ref{mainUP} and Theorem \ref{main2} in
hand, the following proof for  Theorem \ref{voterasy} are exactly
the same with that in \cite{[CP04]}. We left it to the interested
readers. The intuition is that the underlying motion has nothing to
do with the total mass process.

\

\bigskip

\

\bigskip


\end{document}